\documentclass[12pt]{amsart}
\usepackage[english]{babel}
\usepackage[all]{xy}
 \usepackage[latin1]{inputenc}
  \usepackage[T1]{fontenc}
  \usepackage{amsmath,amssymb}
  \usepackage{amsmath}
  \usepackage{amsthm}
\usepackage{amssymb}
\usepackage{mathrsfs}
\usepackage{enumerate}
\usepackage{graphicx}
\usepackage{tikz}\parindent=0.cm
\usetikzlibrary{shapes.geometric}
\usetikzlibrary{positioning}
\usetikzlibrary{arrows}
\usetikzlibrary{shapes.multipart}
\usepackage{physics}
\usepackage[notcite, final, notref]{showkeys}
\usepackage{dsfont}

\newtheorem{theorem}{Theorem}[section]

\newtheorem{lemma}[theorem]{Lemma}
\newtheorem{proposition}[theorem]{Proposition}

\newtheorem{corollary}[theorem]{Corollary}
\newtheorem{definition}[theorem]{Definition}

\theoremstyle{definition}
\newtheorem{remark}[theorem]{Remark}

\newtheorem{example}[theorem]{Example}

\newcommand{\e}{\varepsilon}

\newcommand{\ct}{\circledast}

\usepackage{tikz}
\usetikzlibrary{arrows,shapes,positioning,shadows,trees,matrix}
\usepackage{xcolor}

\begin{document}

\title[Unitary rational functions: The scaled quaternion case]{Unitary rational functions: The scaled quaternion case}

\author[D. Alpay]{Daniel Alpay}
\address{(DA) Schmid College of Science and Technology \\
Chapman University\\
One University Drive
Orange, California 92866\\
USA}
\email{alpay@chapman.edu}

\author[I. Cho]{Ilwoo Cho}
\address{(IC) Department of Mathematics and Statistics \\
Saint Ambrose University \\
508 W. Locust St.
Davenport, IA 52803\\
USA}
\email{choilwoo@sau.edu}

\author[M. Vajiac]{Mihaela Vajiac}
\address{(VA) Schmid College of Science and Technology \\
Chapman University\\
One University Drive
Orange, California 92866\\
USA}
\email{mbvajiac@chapman.edu}

\keywords{hypercomplex, reproducing kernel, quaternions, scaled quaternions, rational functions, realizations}%
\subjclass[2020]{30G35, 46E22} 
\thanks{D. Alpay thanks the Foster G. and Mary McGaw Professorship in
  Mathematical Sciences, which supported his research.}

\maketitle

\begin{abstract}
  We develop the theory of minimal realizations and factorizations of rational functions where the coefficient space is a ring of the type introduced in our
  previous work,  the scaled quaternions, which includes as special cases the quaternions and the split quaternions.
  The methods involved are not a direct generalization of the complex or quaternionic settings, and in particular, the adjoint is not the  classical adjoint
  and we use properties of real Hilbert spaces. This adjoint allows to define the counterpart of unitarity for matrix-rational functions, and we develop the
  corresponding theories of realizations and unitary factorizations. We also begin a theory of matrices in the underlying rings.
\end{abstract}

\tableofcontents

\mbox{}\\

\section{Introduction}
\setcounter{equation}{0}

The theory of linear systems has a long history and we refer for instance to the books \cite{Fuhrmann,MR0255260,MR325201,Zadeh-desoer,
  MR918977, zemanian} for information, background and \footnote{Mathscinet lists almost 400 books with the words {\sl linear systems} 
  in the title.}references.
In view of the associated state space method (see for instance \cite{ot_1_new},\cite{ot21}),
the notion of rational function and the realization theory of matrix-valued rational functions play a key role in linear system theory.
Certain families of matrix-valued rational functions, connected to the notion of dissipativity, are of special interest. To recall the definition
of these families, let $n\in\mathbb N$ and let $J\in\mathbb C^{n\times n}$ be both Hermitian and unitary, that is a signature matrix,
\begin{equation}
  J=J^*=J^{-1}.
  \end{equation}
Furthermore, let $\Omega$ denote either the open unit disk $\mathbb D$, or the right open half-plane $\mathbb C_r$, with
boundary $\partial \Omega$ (thus $\Omega=\mathbb T$, the unit circle, in the first case and $\Omega=i\mathbb R$, the imaginary line, in the
second case). Among the families of rational functions, $\mathbb C^{n\times n}$-valued rational functions $R(z)$, with set of analyticity ${\rm Hol}(R)$,  which satisfy
\begin{eqnarray}
    \label{J-contractive}
    R(z)JR(z)^*  &\le& J, \quad z\in\Omega\cap {\rm Hol}(R)\\
    R(z)JR(z)^*&=&J,\quad z\in\partial\Omega\cap {\rm Hol}(R)
\label{J-unitary-123}
\end{eqnarray}
are called (rational) $J$-inner, and hold a special place. When $J=1$ these are, up to a multiplicative unitary constant, just finite Blaschke factors (respectively of the form 
$\frac{z-a}{z+\overline{a}}$ or $\frac{z-a}{1-z\overline{a}}$ for the $\mathbb C_r$ and $\mathbb D$ cases). In the matrix-valued case, and when
the signature of $J$ is indefinite, $R$ is the product of three kinds of factors, called Blaschke-Potapov factors of the first, second and third kind (the latter are also
called Brune sections). See \cite{pootapov} for the original work, and \cite{ag}, \cite{Dym_CBMS} for later accounts. It is interesting to note that, in the $\mathbb C^{2\times 2}$-case, Brune sections are independently defined and products of Brune considered, in \cite[Probl\`eme 110 p 116]{dbbook}.

\begin{remark}
Condition \eqref{J-unitary-123} can be rewritten as
\begin{eqnarray}
  \label{x-z-R}
  R(z)JR(-\overline{z})^*&=&J,\quad \Omega=\mathbb C_r,\\
  R(z)JR(1/\overline{z})^*&=&J,\quad \Omega=\mathbb T,
                              \label{x-z-t}
\end{eqnarray}

which are equivalent, by analytic continuation, to
\begin{eqnarray}
  \label{x-z-R-1}
  R(x)JR(-x)^*&=&J,\quad \Omega=\mathbb C_r,\\
  R(x)JR(1/x)^*&=&J,\quad \Omega=\mathbb T,
                              \label{x-z-t-1}
\end{eqnarray}
This last definition makes sense in larger settings, and was used in the quaternionic case; see
\cite{zbMATH06658818}. This is the definition we take here for $J$-unitarity in the corresponding ``circle'' and
``line'' cases.
\label{remark-r-t-1}
\end{remark}
Another important family consists of the $\mathbb C^{n\times n}$-valued rational functions $R(z)$ satisfying
\begin{eqnarray}
    \label{J-contractive-2}
    R(z)+R(z)^*&\ge&0  , \quad z\in\Omega\cap {\rm Hol}(R)\\
  R(z)+R(z)^*&=&0,\quad z\in\partial\Omega\cap {\rm Hol}(R).
                \label{rr-2}
\end{eqnarray}

\begin{remark}
One can replace \eqref{rr-2} by
\begin{eqnarray}
  \label{123d}
  R(z)&=&-R(-\overline{z})^*,\quad \Omega=\mathbb C_r,\\
  R(z)&=&-R(1/\overline{z})^*,\quad \Omega=\mathbb T,
          \label{123e}
                     \end{eqnarray}
\label{x-t-r}
and by analytic continuation   \eqref{123d} and \eqref{123e} now take the forms
  \begin{eqnarray}
    \label{phi-1}
    R(x)&=&-R(-x)^*\\
    \label{phi-2}
    R(x)&=&-R(1/x)^{*}.
  \end{eqnarray}
  \label{remark-r-t-2}
\end{remark}

In the setting of complex numbers, condition \eqref{J-contractive} is equivalent to the fact that the kernel
\begin{equation}
K_{R}(z,w)=\frac{J-R(z)JR(w)^*}{1-z\overline{w}}
  \end{equation}
  is positive definite in $\mathbb D\cap{\rm Hol}(R)$ for $\Omega=\mathbb D$ and that the kernel
\begin{equation}
K_{R}(z,w)=\frac{J-R(z)JR(w)^*}{z+\overline{w}}
  \end{equation}
  is positive definite in $\mathbb C_r\cap{\rm Hol}(R)$ for $\Omega=\mathbb C_r$.  Such functions are called {\sl  chain scattering matrices} in network theory, see
e.g. \cite{MR85m:94012}.  In  the setting of complex numbers, condition \eqref{J-contractive} is equivalent to the fact that the kernels
\begin{equation}
L_{R}(z,w)=\frac{R(z)+R(w)^*}{1-z\overline{w}}
  \end{equation}
for $\Omega=\mathbb D$ and 
\begin{equation}
L_{R}(z,w)=\frac{R(z)+R(w)^*}{z+\overline{w}}
  \end{equation}
  for $\Omega=\mathbb C_r$, are positive definite. In network theory, such functions are called {\sl impedance functions}.\smallskip

  As explained in \cite{ag}, these two cases can be reduced to \eqref{J-contractive} with
  \begin{equation}
    \label{JJJJ}
  J=\begin{pmatrix}0&I_{n}\\ I_{n}&0\end{pmatrix}\quad{\rm and}\quad R\quad \mbox{\rm replaced by}\quad \begin{pmatrix}I_{n}&R\\0&I_{n}\end{pmatrix}.
  \end{equation}

  These kernels have also been studied in the case where they have a finite number of negative squares by Krein and Langer in a
  long and fundamental series of papers, of which we mention \cite{kl1,MR47:7504}.\smallskip
  
  As alluded to above, an important tool to study such rational functions is realization theory, which states that any matrix-valued rational
  function with no pole at the origin can be written in the form
  \begin{equation}
    \label{mini-real-1}
R(z)=D+zC(I_N-zA)^{-1}B,
    \end{equation}
    where  $A,B,C$ and $D$ are matrices with complex components and of appropriate sizes. We will also refer to the associated block matrix
    \begin{equation}
      \label{node}
      \begin{pmatrix}A&B\\C&D\end{pmatrix}
      \end{equation}
    as a {\em realization}, or a {\em node}.\smallskip

    Equivalently, if the function is analytic at infinity, the realization is now of the form

    \begin{equation}
    \label{mini-real-2}
R(z)=D+C(zI_N-A)^{-1}B.
\end{equation}

\begin{example}
  Let $w\in\mathbb C$ and $N\in\mathbb N$. Let $R(z)=\frac{1}{(z-w)^N}$. Then, $R$ is of the form \eqref{mini-real-2} with $D=0$ and
\[
  C=\begin{pmatrix}1&0&\cdots &0&0\end{pmatrix}\in\mathbb C^{1\times N},\quad B=\begin{pmatrix}0\\
    0\\ \vdots\\ 0\\1\end{pmatrix}\in\mathbb C^N,
\]
and
\[
  A=\begin{pmatrix}w&1&0&0&0&\cdots\\
    0&w&1&0&0&\cdots\\
    0&0&w&1&0&\cdots\\
    \vdots&&&&&\\
    \vdots&&&&&\\
    0&\cdots &0&0&w&1\\
    0&0&\cdots &0&0&w
  \end{pmatrix}
\]
is a Jordan cell. See e.g. \cite[Exercise 12.2.2 p.  477]{MR3560222}.
\end{example}

      We note that more general type of realizations are possible, but we will not consider them in the present work.
      The McMillan degree of a rational function (say, analytic at the origin) can be defined in a number of equivalent way. It is in particular the minimal  $N$
      in \eqref{mini-real-1} for which a realization exists. The corresponding realization is then called minimal.
      Given a $\mathbb C^{n\times n}$-valued rational function $R$, a factorization $R=R_1R_2$ of $R$ into a product of two
      $\mathbb C^{n\times n}$-valued rational functions $R_1$ and $R_2$, is called {\em minimal} if their degrees add up to the degree of $R$.\smallskip

    Metric properties of the function $R$ translate into metric properties on the node \eqref{node}.
   The following result was proved in \cite[Theorem 2.1 p. 179]{ag},
  see also \cite[Theorem 2.3.1 p. 28]{zbMATH06658818}. In these papers the realization is centered at $\infty$ rather than at the origin. Since the
  imaginary line is invariant under the map $z\mapsto 1/z$, one gets the result by the substitutions $z\mapsto 1/z$ and $w\mapsto 1/w$.
  \begin{theorem}
    Let $J\in\mathbb C^{n\times n}$ be a signature matrix, and let $R(z)$ be a $\mathbb C^{n\times n}$-valued rational function
    regular at the origin, and with realization \eqref{mini-real-1}, $R(z)=D+zC(I_N-zA)^{-1}B$, assumed minimal. Then the following are equivalent:\\
    $(1)$ $R$ takes $J$-unitary-values on the imaginary line, or equivalently
    \begin{equation}
      \label{J-unit}
      R(z)JR(-\overline{z})^*=J,\quad z\in{\rm Hol} (R)\cap {\rm Hol}(\widetilde{R}).
    \end{equation}
    $(2)$ It holds
    \begin{equation}
      DJD^*=J
      \label{lyap-3}
    \end{equation}
    and there exists an Hermitian invertible matrix $H$, uniquely determined from the given realization of $R$, and such that
    \begin{eqnarray}
      \label{lyap}
      A^*H+HA&=&-C^*JC,\\
      C&=&-DJB^*H.
           \label{lypa-2}
      \end{eqnarray}
    \end{theorem}

Then, one has
\begin{equation}
  \label{RJR*}
  \frac{J-R(z)JR(w)^*}{z+\overline{w}}=C(I_N-zA)^{-1}H^{-1}(I_N-wA)^{-*}C^*,
\end{equation}
and in particular $K(z,w)=  \frac{J-R(z)JR(w)^*}{z+\overline{w}}$ is the reproducing kernel of a finite dimensional reproducing kernel Pontryagin space of index $\nu_-(H)$, the number of strictly negative eigenvalues of $H$. Note that the Lyapunov equation \eqref{lyap} need not have a solution, let alone a unique solution. When a solution
exists,  adding the conditions \eqref{lyap-3} and \eqref{lypa-2} ensures uniqueness.\smallskip

One of the main topics considered in \cite{ag} is the concept of minimal unitary factorizations. Given a matrix-valued rational function $R$
satisfying \eqref{J-unit}, one wishes to characterize all, if any, factorizations $R=R_1R_2$ where $R_1$ and $R_2$ have the size as $R$ and also satisfy \eqref{J-unit}
and where the degrees add up.\smallskip

  It is of much interest to replace the field of complex numbers by a general ring, which can be commutative or non-commutative.
  For the general theory, this has been studied in particular by \cite{MR839186,MR0452844,SSR}. Schur analysis in the quaternionic setting
  was considered in a series of papers which includes \cite{MR3568012,acls_milan,MR2872478,acs-survey} and the books \cite{zbMATH06658818,MR4292259}.
  The purpose of the present paper is to study rational functions with some constraints in the setting of
  the family of rings introduced in \cite{MR4591396,MR4592128} and studied in our previous paper \cite{adv_prim}.\smallskip

  In the quaternionic setting, an important tool was the spectral theorem for quaternionic Hermitian matrices. Such a result is not available here,
  and the above definitions have to be adapted in an appropriate way. In the study of rational functions with metric constraints, and as
  for the quaternionic setting, the starting point will be the $\mathbb H_t$-counterpart of \eqref{x-z-R} and of \eqref{x-z-t}.\smallskip
  \begin{remark}
    In the classical setting of the complex numbers, or more generally of a field over the real, the expressions
    \[
(x-a)\quad{\rm and}\quad -(1-xb)
    \]
    are equivalent with $b=a^{-1}$.
An important aspect of linear systems over a ring which admits inverses is that both expressions need to be considered, in particular in the setting of Blaschke factors; see Example \ref{bla-bla}.
    \end{remark}

\begin{remark} This paper is aimed at various audiences (and in particular hypercomplex analysis and classical operator theory), and in consequence we provide sometimes
  more details than necessary.
\end{remark}

The paper consists of six sections besides the introduction.  In Section \ref{2} we briefly survey the main aspects of the rings $\mathbb H_t$.  Realization theory of 
rational functions with coefficients in $\mathbb H_t$ is developed in Section \ref{sec:real}. The cases of rational functions with metric constraints are considered in the 
remaining sections. The table of contents provides a more detailed picture, and we direct the reader to the title of the subsections for clarity.

\section{The rings $\mathbb H_t$}
\setcounter{equation}{0}
\label{2}
\subsection{Definition and some lemmas}
In \cite{MR4591396,MR4592128,alpay_cho_3} the following algebras were introduced and studied.
\begin{definition}
  Let $t\in\mathbb R$. We denote by $\mathcal H_2^t$ the set of matrices of the form
  \[
    \begin{pmatrix}a&tb\\ \overline{b}&\overline{a}\end{pmatrix},
  \]
  with $a,b\in\mathbb C$.
\end{definition}

For the following result we refer to \cite{MR4591396,MR4592128,alpay_cho_3}.

\begin{theorem}
  $\mathcal H_2^t$ is an algebra for all $t\in\mathbb R$ when endowed with the usual matrix operations. It is a skew field for $t<0$ and in particular
  invariant under inversion. When $t>0$, any element of $\mathcal H_2^t$ invertible in $\mathbb C^{2\times 2}$ has its inverse in $\mathcal H_2^t$.
  \end{theorem}

When $t=-1$ we obtain the quaternions, while $t=1$ corresponds to the split quaternions.
We note that $\mathcal H_2^t$ is not closed under the usual matrix adjoint for $t\not=\pm 1$, and this is an important difference with the quaternions and the split
quaternions, and require different and more general methods.

\begin{remark}
  As pointed out to us by Dan Volok, one has for $t>0$,
  \begin{equation}
    \begin{pmatrix}a&tb\\ \overline{b}&\overline{a}\end{pmatrix}=   \begin{pmatrix}1&0\\ 0& \sqrt{t}^{-1}\end{pmatrix}
         \begin{pmatrix}a& \sqrt{t}b\\ \sqrt{t}\,\overline{b}&\overline{a}\end{pmatrix}  \begin{pmatrix}1&0\\ 0& \sqrt{t}\end{pmatrix},
       \end{equation}
       while, for $t<0$,
  \begin{equation}
    \begin{pmatrix}a&tb\\ \overline{b}&\overline{a}\end{pmatrix}=   \begin{pmatrix}1&0\\ 0&\sqrt{-t}^{-1}\end{pmatrix}
         \begin{pmatrix}a& -\sqrt{-t}b\\ \sqrt{-t}\,\overline{b}&\overline{a}\end{pmatrix}  \begin{pmatrix}1&0\\ 0& \sqrt{-t}\end{pmatrix}.
       \end{equation}
  \end{remark}

The following lemma holds trivially for $t<0$ since then $\mathcal H^t_2$ is a skew-field.
\begin{lemma}
  \label{debut0}Let $t\not=0$ and assume that
  \begin{equation}
    \label{debut}
\begin{pmatrix}a&tb\\ \overline{b}&\overline{a}\end{pmatrix}\begin{pmatrix}u\\ \overline{v}\end{pmatrix}=\begin{pmatrix}0\\ 0\end{pmatrix}.
\end{equation}
Then
\begin{equation}
  \label{fin}
  \begin{pmatrix}a&tb\\ \overline{b}&\overline{a}\end{pmatrix}
  \begin{pmatrix}tv\\ \overline{u}\end{pmatrix}=\begin{pmatrix}0\\ 0\end{pmatrix}.
\end{equation}
Similarly, assume that
\begin{equation}
  \label{debut2}
  \begin{pmatrix}u&tv\end{pmatrix}\begin{pmatrix}a&tb\\ \overline{b}&\overline{a}\end{pmatrix}=\begin{pmatrix}0& 0\end{pmatrix}.
\end{equation}
Then,
\begin{equation}
  \label{fin2}
\begin{pmatrix}\overline{v}&\overline{u}\end{pmatrix}\begin{pmatrix}a&tb\\ \overline{b}&\overline{a}\end{pmatrix}=\begin{pmatrix}0& 0\end{pmatrix}.
\end{equation}
  \end{lemma}
  \begin{proof}
    \eqref{debut} is equivalent to
    \[
      \begin{split}
        au+tb\overline{v}&=0,\\
        \overline{b}u+\overline{a}\overline{v}&=0,
        \end{split}
      \]
      while \eqref{fin} is equivalent to
          \[
      \begin{split}
        atv+tb\overline{u}&=0,\\
       \overline{b}tv+\overline{a}\overline{u}&=0,
        \end{split}
      \]
      which are equivalent since $t\not=0$.\smallskip

      In a similar way, \eqref{debut2} is equivalent to
      \begin{equation}
        \label{fin3}
        \begin{split}
          ua+tv\overline{b}&=0,\\
          utb+tv\overline{a}&=0,
        \end{split}
      \end{equation}
      and \eqref{fin2} is equivalent to
      \begin{equation}
        \begin{split}
        \overline{v}a+\overline{u}\overline{b}&=0,\\
        \overline{v}tb+\overline{u}\overline{a}&=0,
      \end{split}
      \label{fin4}
    \end{equation}
    which is equivalent to \eqref{fin3} since $t\not=0$.
  \end{proof}
Rather than a matrix notation we will use a notation representing $\mathcal H_2^t$ as a vector space over the real numbers. Writing
$a=x_0+ix_1$ and $b=x_2+ix_3$ (where $x_0,x_1,x_2$ and $x_4$ are real numbers) and
\begin{equation}
  \label{real-matrix-form}
  \left(\begin{array}{cc}a & tb\\
\overline{b} & \overline{a}
        \end{array}\right)
=x_0\left(\begin{array}{cc}
1 & 0\\
0 & 1
\end{array}\right)+
x_1
\left(\begin{array}{cc}
i & 0\\
0 & -i
\end{array}\right)+
x_2\left(\begin{array}{cc}
0 & t\\
1 & 0
\end{array}\right)+
x_3\left(\begin{array}{cc}
0 & ti\\
-i & 0
\end{array}\right),
\end{equation}
we set
\begin{equation}
  i_t=\left(\begin{array}{cc}
i & 0\\
0 & -i
             \end{array}\right),\quad   j_t=\left(\begin{array}{cc}
0 & t\\
1 & 0
             \end{array}\right),\quad            k_t=\left(\begin{array}{cc}
0 & ti\\
-i & 0
\end{array}\right)
\end{equation}
with Cayley table of multiplication

\begin{center}
  \begin{tabular}{|c|c|c|c|c|}
    \hline
    $\nearrow$&$1$&$i_t$& $j_t$&$k_t$\\
    \hline
    $1$&$1$&$i_t$& $j_t$&$k_t$\\
    \hline
    $i_t$&$i_t$&$-1$ &$k_t$ &$-j_t$\\
    \hline
    $j_t$& $j_t$&$-k_t$&$t$&$-ti_t$ \\
    \hline
    $k_t$&$k_t$ &$j_t$ &$ti_t$ &$t$\\
    \hline
  \end{tabular}
\end{center}

\begin{definition}
We denote
\begin{equation}
  \label{ht}
  \mathbb H_t=\left\{x_0+x_1i_t+x_2j_t+x_3k_t,\quad x_0,x_1,x_2,x_3\in\mathbb R\right\}.
    \end{equation}
\end{definition}


\begin{remark}
  We write $i$ rather than $i_t$,   and will remove the dependence on $t$ in the product, meaning that we view the ring $\mathbb H_t$ as containing a fixed copy of the complex numbers. The above Cayley table and \eqref{ht}
  thus become
\begin{center}
  \begin{tabular}{|c|c|c|c|c|}
    \hline
    $\nearrow$&$1$&$i$& $j_t$&$k_t$\\
    \hline
    $1$&$1$&$i_t$& $j_t$&$k_t$\\
    \hline
    $i$&$i$&$-1$ &$k_t$ &$-j_t$\\
    \hline
    $j_t$& $j_t$&$-k_t$&$t$&$-ti$ \\
    \hline
    $k_t$&$k_t$ &$j_t$ &$ti$ &$t$\\
    \hline
  \end{tabular}
\end{center}

and
\begin{equation}
  \label{ht1}
  \mathbb H_t=\left\{x_0+x_1i+x_2j_t+x_3k_t,\quad x_0,x_1,x_2,x_3\in\mathbb R\right\}.
    \end{equation}
\end{remark}

We use the notation
\begin{equation}
  \label{Iq}
  I(q)=\begin{pmatrix}a&tb\\ \overline{b}&\overline{a}\end{pmatrix},\quad q=a+bj_t\in\mathbb H_t,
\end{equation}
and for $A=(a_{jk})\in\mathbb H_t^{n\times m}$
\begin{equation}
  \label{I-A}
I(A)=(I(a_{jk})).
\end{equation}

\begin{lemma}
  \label{lemm-ab-I}
  Let $A\in\mathbb H_t^{n\times m}$ and $B\in\mathbb H_t^{m\times r}$. It holds that
  \begin{equation}
    I(AB)=I(A)I(B).
    \end{equation}
  \end{lemma}

  \begin{proof}
    The result is true for $n=m=1$.  Let now $A$ and $B$ as in the statement of the lemma. The matrices $I(A)$ and $I(B)$ belong to $\mathbb C^{2n\times 2m}$ and $\mathbb C^{2m\times 2r}$ respectively. Since multiplication in block by matrices does not depend on compatible block decompositions
    (see e.g. \cite[Theorem 1.3.2 p. 27]{MR1002570}), we have:
    \[
      \begin{split}
        I(A)I(B)&=(I(a_{uv}))(I(b_{vw}))\\
        &=(\sum_{v=1}^mI(a_{uv})I(b_{vw}))\\
        &=(\sum_{v=1}^mI(a_{uv}b_{vw}))\\
        &=(I(c_{uw})),\qquad{\rm with}\quad c_{uw}=\sum_{v=1}^ma_{uv}b_{vw}\\
        &=I(C)
      \end{split}
    \]
    where we have used the fact that the claim is true for scalars, to go from the third line to the fourth line. 
    \end{proof}

    As a direct corollary of the previous lemma we have:
    \begin{corollary}
      \label{ab-ba}
    Let $A\in\mathbb H_t^{n\times m}$ and $B\in\mathbb H_t^{m\times n}$ be such that
    \begin{equation}
\label{prod-ab-ba}
      AB=I_{\mathbb H_t^n}\quad and\quad      BA=I_{\mathbb H_t^m}.
    \end{equation}
    Then, $n=m$, and (by definition), $A=B^{-1}$ and $B=A^{-1}$.
  \end{corollary}

  \begin{proof} We have that $I(A)\in\mathbb C^{2n\times 2m}$ and $I(B)\in\mathbb C^{2m\times 2n}$.
Using Lemma \ref{lemm-ab-I} we see that \eqref{prod-ab-ba} is equivalent to
    \[
I(A)I(B)=I_{2n}\quad {\rm and}\quad I(B)I(A)=I_{2m}.
\]
from which $2n=2m$ since  $I(A)$ and $I(B)$ are matrices with complex entries and so $n=m$.
\end{proof}

As a consequence of the above we have
\begin{equation}
I(A^{-1})=(I(A))^{-1},
\end{equation}
meaning that if one of the sides of the equality makes sense so does the other, and then both sides coincide.

\begin{proposition} Let $t\not=0$ and  $M,N\in\mathbb N$.
  The space $\mathbb H_t^N$ is a $\mathbb H_t$ right and left module and has a basis. Thus every right (or left) linear map from $\mathbb H_t^N$ into
  $\mathbb H_t^M$  can be represented by a matrix.
\end{proposition}

\begin{proof}
 See \cite[A II. 143]{MR0274237}.
\end{proof}

\begin{proposition}
  \label{prop-eigen}
  Let $A\in\mathbb H_t^{n\times n}$. There exist $f\in\mathbb H_t^{n}$ and $\lambda\in \mathbb H_t$ such that
  \begin{equation}
    \label{2-21}
Af=f\lambda.
    \end{equation}
  \end{proposition}

  \begin{proof}
    The matrix $I(A)\in\mathbb C^{2n\times 2n}$, and so has, by the fundamental theorem of algebra, at least one complex eigenvalue, say $\lambda$, with associated eigenvector $u\in\mathbb C^{2n}$, that is
    \begin{equation}
      I(A)u=\lambda u.
      \label{chi-A-u}
      \end{equation}
    We write $u$ in the form
    \[
      u=\begin{pmatrix} a_1\\ \overline{b_1}\\ a_2\\ \overline{b_2}\\ \vdots\\ a_n\\ \overline{b_n}\end{pmatrix},\quad\mbox{\rm and define}\quad v
=\begin{pmatrix} tb_1\\ \overline{a_1}\\ tb_2\\ \overline{a_2}\\ \vdots\\ tb_n\\ \overline{a_n}\end{pmatrix}.
\]
Note that $v$ is constructed in such a way that there exists (a uniquely defined) $f\in\mathbb H_t^n$ such that
$I(f)=\begin{pmatrix}u&v\end{pmatrix}$.
With $A=(A_{ij})$ and $I(A_{ij})=\begin{pmatrix} \alpha_{ij}&t\beta_{ij}\\ \overline{\beta_{ij}}&\overline{\alpha_{ij}}\end{pmatrix}$, \eqref{chi-A-u} can be rewritten as
\[
\sum_{j=1}^n\begin{pmatrix} \alpha_{ij}&t\beta_{ij}\\ \overline{\beta_{ij}}&\overline{\alpha_{ij}}\end{pmatrix}\begin{pmatrix}a_j\\ \overline{b_j}\end{pmatrix}=\lambda \begin{pmatrix}a_i\\ \overline{b_i}\end{pmatrix}, \quad i=1,\ldots, n,
\]
that is
\begin{eqnarray}
  \label{eq1234}
  \sum_{j=1}^n\alpha_{ij} a_j+t\beta_{ij}\overline{b_j}&=&\lambda a_i,\\
  \label{eq12345}
  \sum_{j=1}^n\overline{\beta_{ij}} a_j+ \overline{\alpha_{ij}}\overline{b_j}&=&\lambda \overline{b_i},                                                                                                                        \quad i=1,\ldots, n.
  \end{eqnarray}
Using \eqref{eq1234}-\eqref{eq12345} we now show that
\begin{equation}
  \label{chi-A-v}
  Av=\overline{\lambda}v.
\end{equation}
Indeed, \eqref{chi-A-v} can be rewritten as
\[
  \sum_{j=1}^n\begin{pmatrix} \alpha_{ij}&t\beta_{ij}\\ \overline{\beta_{ij}}&\overline{\alpha_{ij}}\end{pmatrix}\begin{pmatrix}tb_j\\ \overline{a_j}
  \end{pmatrix}=\overline{\lambda} \begin{pmatrix}b_i\\ \overline{a_i}\end{pmatrix}, \quad i=1,\ldots, n,
\]
that is
\[
  \begin{split}
    \sum_{j=1}^n \alpha_{ij}b_j+\beta_{ij} \overline{b_j}&=\overline{\lambda}{b_i},\\
    \sum_{j=1}^n\overline{\beta_{ij}}b_j+\overline{\alpha_{ij}}\overline{ a_j}&=\overline{\lambda}\overline{ a_i},   \quad i=1,\ldots, n.
      \end{split}
  \]
  but these are just the conjugates of \eqref{eq12345} and \eqref{eq1234} respectively.
  Combining \eqref{chi-A-u} and \eqref{chi-A-v} we obtain
  \[
I(A)\begin{pmatrix}u&v\end{pmatrix}=\begin{pmatrix}u&v\end{pmatrix}\begin{pmatrix}\lambda&0\\0&\overline{\lambda}\end{pmatrix},
\]
which allows to conclude with $f\in\mathbb H_t^n$ defined as above by $I(f)=\begin{pmatrix}u&v\end{pmatrix}$ since $I(\lambda)=
\begin{pmatrix}\lambda&0\\0&\overline{\lambda}\end{pmatrix}$ for $\lambda\in\mathbb C$.
\end{proof}

When $t>0$, the existence of divisors of $0$ prevents from calling $f$ an eigenvector and $\lambda$ an eigenvalue of $A$.  We note that the one dimensional real vector space spanned by $f$ {\sl is not} $A$-invariant, but in special cases (for instance if $\lambda$ is real or if $f\lambda=0$, the latter being
possible only for $t>0$).

\begin{lemma} Let $f$ and $\lambda$ be as in \eqref{2-21}.
  The space
  \begin{equation}
    \mathfrak M_f=\left\{f\lambda^n\quad\vert\,\, n=0,1,2,\ldots\right\}
  \end{equation}
  is a finite dimensional real vector space which is $A$-invariant. It is of dimension $1$ if $\lambda\in\mathbb R$ and of dimension $2$
  if $\lambda\not=\overline{\lambda}$.
    \end{lemma}

  \begin{proof}
    This comes from the associativity of multiplication of matrices (of possibly different sizes) with entries in $\mathbb H_t$,
    \[
A^2f=A(Af)=A(f\lambda)=(Af)\lambda=f\lambda^2,\quad A^3f=f\lambda^3,\ldots
\]
Furthermore, the space is at most two dimensional since the real linear span of $1,\lambda,\ldots $ is at most two dimensional. This follows from the fact that the roots
of the equation
\[
z^2-2({\rm Re}\,\lambda) z+|\lambda|^2=0
\]
are $z=\lambda$ and $z=\overline{\lambda}$. The claim on the dimension of $\mathfrak M_f$ follows,
    \end{proof}

When $p\not=\pm 1$, $\mathbb H_t$ is not invariant under matrix-conjugation, but two adjoints which leave $\mathcal H_2^t$ invariant can be defined,
namely $q^\ct$ and $q^{[*]}$, defined by
\begin{equation}
  \label{q1}
I(  q^{\circledast})=\begin{pmatrix}\overline{a}&-tb\\ -\overline{b}&a\end{pmatrix}
\end{equation}
and
\begin{equation}
  \label{q2}
I(  q^{[*]})=\begin{pmatrix}a&t\overline{b}\\ b&\overline{a}\end{pmatrix}.
  \end{equation}

  In the present paper we focus on the first adjoint and leave the parallel study for the second adjoint to a future publication.

  \subsection{Point evaluations}

  Given a sequences $(f_n)_{n\in\mathbb N_0}$ of elements of $\mathbb H_t$ one can define the power series
    \[
f(x)=\sum_{n=0}^\infty f_nx^n,
    \]
    which may, or may not, converge in a neighborhood of the origin. The point-wise product of $f$ with another such power series, say
    $g(x)=\sum_{n=0}^\infty g_nx^n$, is not commutative, and correspond to the
    convolution, also called Cauchy product (see e.g. \cite{MR51:583}), of the sequences $(f_n)_{n\in\mathbb N_0}$ and $(g_n)_{n\in\mathbb N_0}$:
    \[
h_n=\sum_{k=0}^n f_kg_{n-k},\quad n=0,1,\ldots
\]
\begin{definition}
  For $f(x)=\sum_{n=0}^\infty f_nx^n$ as above and $q\in\mathbb H_t$ we get the point evaluation on the right
  \begin{equation}
    f(q)=\sum_{n=0}^\infty f_nq^n.
  \end{equation}
  For $f$ and $g$ as above, we define the $\star$ product to be
  \begin{equation}
(f\star g)(q)=\sum_{n=0}^\infty\left(\sum_{k=0}^n f_kg_{n-k}\right)q^n.
    \end{equation}
  \end{definition}
\subsection{The adjoint $\ct$}
    As mentioned above, there are two natural adjoints in the algebra $\mathbb H_t$, denoted in our previous paper \cite{acv1} by $\ct$ and $[*]$; see \eqref{q1} and \eqref{q2}. In this paper we focus on $\ct$, defined by
   (see \cite[(2.16)]{adv_prim} and \eqref{q1} above)
    \[
p^\ct=\overline{a}-j_tb,
\]
with associated real bilinear form
\begin{equation}
  \label{inner-prod-h-t}
[p,q]_{\circledast}={\rm Tr}\, (I(p^{\circledast})I(q))=a\overline{c}+\overline{a}c-t(b\overline{d}+\overline{b}d),
      \end{equation}
      where $q=a+j_tb$ and $p=c+j_td$. It is in particular real-valued. A very important feature of $\ct$, not shared by $[*]$, is the following fact, stated
      as a lemma. The proof is a direct computation and will be omitted. See also \cite{acv1}.

  The proofs of the following lemmas are simple and will be omitted.

  \begin{lemma}
    \label{unit-p-0}
     Let $p=a+j_t b\in\mathbb H_t$. Then,
        \[
          pp^\ct=|a|^2-t|b|^2
        \]
        is real. Furthermore, 
    $pp^\ct=1$ if and only if
    \begin{equation}
|a|^2-t|b|^2=1.
      \end{equation}
    \end{lemma}

     \begin{lemma}
       \label{unit-p-0-0}
            Let $p=a+j_t b\in\mathbb H_t$. Then,
    $p+p^\ct=0$ if and only if $a$ is purely imaginary.
       \end{lemma}
  
      We note that
      \begin{equation}
        \label{sym-345}
[p,q]_{\circledast}=[q,p]_{\circledast},\quad \forall p,q\in\mathbb H_t.
        \end{equation}

      Although already mentioned, we recall that $t\not=0$. The following lemma fails to hold for $t=0$.

        \begin{lemma}
        Let $p\in\mathbb H_t$ be such that
        \begin{equation}
          [q,p]_\ct=0,\quad \forall p\in\mathbb H_t.
        \end{equation}
        Then, $q=0$.
      \end{lemma}

      \begin{proof}
Let first $d=0$ and $c=1$ to get ${\rm Re}\, a=0$ in \eqref{inner-prod-h-t}. Still for $d=0$ the choice $c=i$ gives ${\rm Im}\, a=0$ and hence $a=0$. The proof that $b=0$ is done in a similar way.  
        \end{proof}

        \begin{lemma}
        \label{lem-1} Let $p,q,h\in\mathbb H_t$. It holds that
  \begin{equation}
    [q,hp]_{\circledast}=[h^{\circledast}q,p]_{\circledast}.
    \end{equation}
  \end{lemma}

  \begin{proof}
    We have
    \[
\begin{split}
  [q,hp]_{\circledast}&={\rm Tr}\,(I((hp)^{\circledast})I(q))\\
  &={\rm Tr}\,(I(p^{\circledast}h^\ct)I(q))\\
  &={\rm Tr}\,(I(p^{\circledast})I(h^\ct)I(q))\\
  &={\rm Tr}\,(I(p^{\circledast})I(h^\ct q))\\
  &={\rm Tr}\,(I(h^\ct q)I(p^{\circledast})  )\\
  &=[h^{\circledast}q,p]_{\circledast}.
\end{split}
\]
\end{proof}

For
\[
  h=\begin{pmatrix}h_1\\ \vdots \\ h_n\end{pmatrix}\quad{\rm and}\quad g=\begin{pmatrix}g_1\\ \vdots\\ g_n\end{pmatrix}\in\mathbb H_t^n,
\]
  we define
  \begin{equation}
    \label{finite-pi-k}
    [h,g]_{\circledast} =\sum_{k=1}^n   [h_k,g_k]_{\circledast}.
  \end{equation}

  \begin{proposition}
    \label{pontry}
    For $t<0$, $\mathbb H_t^n$ endowed with \eqref{finite-pi-k} is a finite dimensional Hilbert space on the real numbers.
For $t>0$, $\mathbb H_t^n$ endowed with \eqref{finite-pi-k} is a finite dimensional Pontryagin space on the real numbers.
    \end{proposition}
    \begin{proof}
      For $t>0$, consider the map $J$ from $\mathcal H_t$ into itself defined by
      \begin{equation}
        \label{J-pont}
        Jp=a-j_tb,\quad p=a+j_tb.
      \end{equation}
      Then, replacing $d$ by $-d$ in \eqref{inner-prod-h-t}, we obtain
      \[
        \begin{split}
          \langle p ,q\rangle_\ct&\stackrel{\rm def}{=}[p,Jq]_\ct\\&=a\overline{c}+\overline{a}c+t(b\overline{d}+\overline{b}d)\\
          &=\frac{1}{2}\left((a+\overline{a})(\overline{c+\overline{c}})+(a-\overline{a})(\overline{c-\overline{c}})+t(b+\overline{b})(\overline{d+\overline{d}})+t(b-\overline{b})(\overline{d-\overline{d}})\right),
          \end{split}
\]
 which    makes $\mathbb H_t$ into a Hilbert space. The claim then holds since $\mathbb H_t^n$ is a finite direct sum of finite dimensional Pontryagin spaces.
\end{proof}

\begin{remark}
  \label{new-remark}
 The matrices $I(q)$ and $I(q^\ct)$ are similar since they have the same spectrum, but the similarity depends on the given point $q$.
  \end{remark}

\begin{definition}
  \label{anti-not-anti}
For $A=(a_{jk})\in\mathbb H_t^{u\times v}$ we define
\begin{equation}
  \label{labek-a-ct}
  A^\ct=(a_{kj}^\ct)\in\mathbb H_t^{v\times u}.
\end{equation}
When $u=v$ the matrix $A$ is called $\ct$-symmetric if $A=A^\ct$ and $\ct$-anti-symmetric if $A=-A^\ct$.
\end{definition}

The decomposition
\begin{equation}
  A=\frac{A+A^\ct}{2}+\frac{A-A^\ct}{2}
\end{equation}
of $A$ into a sum of a $\ct$-symmetric and of a $\ct$-anti-symmetric matrix is the counterpart of the classical case (where one can divide
by $i$ and get two Hermitian matrices) and will be used in what follows.

\begin{remark}
  \label{remstar}
  We note that $I(A^\ct)\not= I(A)^*$ when $t\not=\pm 1$. So in opposition to the quaternionic setting, one cannot reduce directly
  results involving $\ct$ to the complex setting. It is also worth noting that although $I(q)$ and $I(q^\ct)$ are similar for $q\in\mathbb H_t$ one cannot infer a similarity for the matrices $I(A)$ and $I(A^\ct)$  since the similarity for each entry depends on the entry itself.
  \end{remark}

  \begin{lemma}
    $A\in\mathbb H^{n\times n}$ is invertible if and only if $A^\ct$ is invertible, and we have
    \begin{equation}
      (A^{-1})^\ct=(A^\ct)^{-1}.
    \end{equation}
  \end{lemma}
  \begin{proof}
    Let $B\in\mathbb H_t^{n\times n}$ be such that
    \begin{equation}
      \label{abba123}
      AB=I_{\mathbb H_t^n}.
    \end{equation}
    Thus $I(A)I(B)=I_{2n}$ so that $I(B)I(A)=I_{2n}$ and $BA=I_{\mathbb H_t^n}$.
    Applying $\ct$ to \eqref{abba123} we have
    \[
      B^\ct A^\ct=I_{\mathbb H_t^n},
    \]
    and as above $ A^\ct B^\ct=I_{\mathbb H_t^n}$.
    Hence, $(A^\ct)^{-1}=B^\ct=(A^{-1})^\ct$.    
    \end{proof}
  We denote
  \begin{equation}
    \label{A-ct}
    A^{-\ct}:=  (A^{-1})^\ct=(A^\ct)^{-1}.
    \end{equation}
\begin{proposition}
        \label{adj-a}
    Let $A\in\mathbb H_t^{n\times m}$, $q\in\mathbb H_t^m$ and $p\in\mathbb H_t^n$.
    Then,
    \begin{equation}
      [Aq,p]_{\circledast}=      [q,A^{\circledast}p]_{\circledast},
    \end{equation}
    where $A^\ct$ has been defined in \eqref{labek-a-ct}.
    \end{proposition}

    \begin{proof}
      We have
      \[
        (Aq)_\ell=\sum_{j=1}^ma_{\ell j}q_j,\quad \ell=1,\ldots, n,
      \]
      and so
      \[
        \begin{split}
          [Aq,p]_{\circledast}&=\sum_{\ell=1}^n[\sum_{j=1}^ma_{\ell j}q_j,p_\ell]_{\circledast}\\
          &=\sum_{j=1}^m\sum_{\ell=1}^n[a_{\ell j}q_j,p_\ell]_{\circledast}\\
          &=\sum_{j=1}^m[\sum_{\ell=1}^na_{\ell j}q_j,p_\ell]_{\circledast}\\
          &=\sum_{j=1}^m[\sum_{\ell=1}^nq_j,a_{\ell j}^{\circledast}p_\ell]_{\circledast}\quad\mbox{(by Lemma \ref{lem-1})}\\
                    &=\sum_{j=1}^m[q_j,\sum_{\ell=1}^na_{\ell j}^{\circledast}p_\ell]_{\circledast}\\
          &=\sum_{j=1}^m[q_j,(A^{\circledast}p)_j]_{\circledast}\\
          &=          [q,A^{\circledast}p]_{\circledast}.
        \end{split}
        \]
      \end{proof}

      \subsection{Convergence in $\mathbb H_t^n$}
      \label{conv-123}
      With the map $I$ defined by \eqref{Iq} and \eqref{I-A},  we recall (see \eqref{inner-prod-h-t}) that
      \[
[f,g]_{\ct}={\rm Tr}\, \left(\left(I(g^\ct)\right)^tI(f)\right),\quad f,g\in\mathbb H_t^n.
\]
When $t<0$ we have an inner product and a norm and in this case, the convergence is done using the bracket $[\cdot,\cdot]_\ct$ above. When $t>0$ we have a Pontryagin space and can use the
following result of \cite[Theorem 2.4 p. 18]{ikl} in the complex setting case. The result remains true for a real Pontryagin space and is omitted.

\begin{proposition}
  Let $(\mathcal P,[\cdot,\cdot]_{\mathcal P})$ be a Pontryagin space over the real numbers. A sequence $(f_n)_{n\in\mathbb N}$ of elements in $\mathcal P$ converges to
  $f_0\in\mathcal P$ if and only if the following two conditions hold:\\
  \begin{eqnarray}
    \lim_{n\rightarrow\infty}[f_n,f_n]_{\mathcal P}&=&[f_0,f_0]_{\mathcal P}\\
    \lim_{n\rightarrow\infty}[f_n,g]_{\mathcal P}&=&[f_0,g]_{\mathcal P}
                                                     \label{cv-weak}
\end{eqnarray}
for $g$ in a dense set of $\mathcal P$.
  \end{proposition}
  When $\mathcal P$ is a reproducing kernel Pontryagin space, one can take in the second condition the reproducing kernels, and then \eqref{cv-weak} amounts to
  point-wise convergence. This makes the result of interest even in the finite dimensional setting. Using the map $I$ one reduces then the proof of the $\mathbb H_t$ setting to the real setting.
  
\subsection{$\ct$-symmetric matrices}
Let $q=a+bj_t\in\mathbb H_t$ with $a,b\in\mathbb C$. By definition of $\circledast$ we have that $q=q^{\circledast}$ if and only if
\[
a+bj_t=\overline{a}-bj_t,
\]
i.e. if and only if $a\in\mathbb R$ and $b=0$, i.e. if and only if $q$ is a real number.

\begin{definition}
  The matrix $M\in\mathbb H_t^{n\times n}$ is called $\ct$-symmetric if $M=M^\circledast$, and semi-definite $\ct$-positive
  (or $\ct$-non-negative) if
  \begin{equation}
    [p,Mp]_\ct\ge 0,\quad\forall p\in\mathbb H_t^n.
  \end{equation}
  We then will use the notation
  \[
M\ge_\ct 0.
    \]
\end{definition}

In this subsection we focus on the case of $\ct$-symmetric matrices; the case of $\ct$-positive matrices is considered in the next section.
As already mentioned, the space $\mathbb H_t^n$ endowed with the symmetric real-valued bilinear form
\begin{equation}
[f,g]_\ct=\sum_{a=1}^n[f_a,g_a]_\ct,\quad f=\begin{pmatrix}f_1\\ f_2\\ \vdots\\ f_n\end{pmatrix}\quad {\rm and}\quad g=\begin{pmatrix}g_1\\ g_2\\ \vdots\\ g_n\end{pmatrix}
  \end{equation}
  is a Pontryagin space on the real numbers. With $H=H^\ct\in\mathbb H_t^{n\times n}$ invertible, we define a new inner product on $(\mathbb H_t^n,[\cdot,\cdot]_\ct)$ via
  \begin{equation}
    \label{[]H}
    [f,g]_H=[f,Hg]_\ct,\quad f,g\in\mathbb H_t^n.
  \end{equation}

  We denote by $A^\Delta$ the adjoint of a matrix $A\in\mathbb H_t^{n\times n}$ with respect to this form, meaning
  \begin{equation}
[Af,g]_H=[f,A^\Delta g]_H.
\end{equation}

\begin{proposition} It holds that
  \begin{equation}
    \label{Adelta}
A^\Delta=H^{-1}A^\ct H.
    \end{equation}
  \end{proposition}

  \begin{proof} We have
    \[
      \begin{split}
        [Af,g]_H&=[Af,Hg]_\ct\\
        &=[f,A^\ct Hg]_\ct\\
        &=[f,HH^{-1}A^\ct Hg]_\ct\\
        &=[f,H^{-1}A^\ct H]_H
      \end{split}
    \]
    and hence $A^\Delta=H^{-1}A^\ct H$ holds.
  \end{proof}

  It is readily seen that
  \begin{equation}
    (AB)^\Delta=B^\Delta A^\Delta,\quad A,B\in\mathbb H_t^{n\times n}.
  \end{equation}

  For $\mathcal M\subset\mathbb H_t^n$ we define
  \begin{equation}
\mathcal M^{\perp_H}=\left\{f\in\mathbb H_t^n\,\,:\,\, [f,m]_H=0,\quad\forall m\in\mathcal M\right\}.
\end{equation}

  The following result will play a key role in Sections \ref{c45} and \ref{c57}, where minimal factorizations with metric constraints are considered.

  \begin{lemma}
    \label{lemma-preumss}
    In the previous notation, assume that $A\mathcal M\subset\mathcal M$.  Then $A^\Delta\mathcal M^{\perp_H}\subset \mathcal M^{\perp_H}$.
    \end{lemma}
  
  \begin{proof}
    We have for $f\in\mathcal M^{\perp_H}$ and $m\in \mathcal M$,
    \[
      [A^\Delta f,n]_H=[f,\underbrace{An}_{\substack{\in\mathcal M\\(A\mathcal M\subset\mathcal M)}}]_H=0.
      \]
    \end{proof}
The space $\mathcal M$ is non-degenerate with respect to $[\cdot,\cdot]_H$ if $\mathcal M\cap\mathcal M^{\perp_H}=\left\{0\right\}$.
Then we have the direct and orthogonal sum
\begin{equation}
\mathbb H_t^n=\mathcal M\oplus_H\mathcal M^{\perp_H}.
  \end{equation}

\begin{lemma}
There exists $f\in\mathbb H_t^n$ such that $[f,f]_H\not=0$.
\end{lemma}

\begin{proof}
Assume that $[f,f]_H=0$ for all $f\in\mathbb H_t^n$, and let $g\in\mathbb H_t^n$. Then we also have
\[
[f+g,f+g]_H=[g,g]_H=0,
\]
and hence $[f,g]_H=0$ since 
\[
[f+g,f+g]_H=[f,f]_H+[g,g]_H+2[f,g]_H.
\]
Then, $[f,Hg]_\ct=0$ for all $f,g\in\mathbb H_t^n$, which implies $H=0$.
But this is impossible since $H$ is invertible and in particular, $H\not=0$.
\end{proof}

We note that the above result could also be proved using the map $I$ and the polarization formula for complex matrices.

\begin{remark}
By the preceding lemma there exists $f\in\mathbb H_t^n$ such that $[f,f]_H\not=0$. Let
$f_0=\frac{f}{\sqrt{|[f,f]_H|}}$.
We write
\[
  \mathbb H_t^n=\mathcal M_0\oplus_H \mathcal M_0^{\perp_H}.
\]  
\end{remark}

  \subsection{$\ct$-Positive matrices}
  \begin{proposition}
    Let $M\in \mathbb H_t^{n\times n}$ be $\ct$-symmetric. Then,
$M\ge_\ct0$ if and only if it factorizes through a real Hilbert space.
  \end{proposition}

\begin{proof}
  We just outline the proof, which follows classical arguments; see e.g. \cite[Theorem 2.2 p.703]{atv1} for a related result. We consider the range of $M$ endowed with the form
  \begin{equation}
    \label{MpMq}
[Mp,Mq]_M=[p,Mq]_\ct=[Mp,q]_\ct,\quad p,q\in\mathbb H_t^n.
\end{equation}

STEP 1:{\sl The symmetric form  \eqref{MpMq} is well defined and $({\rm ran}\, M,[\cdot,\cdot]_M)$ is a Hilbert space.}\\

Indeed, let $Mp_1=Mp_2$  and  $Mq_1=Mq_2$ for $p_1,p_2,q_1,q_2\in\mathbb H_t^n$. It holds that
\[
  \begin{split}
    [Mp_1,Mq_1]_M&=[p_1,Mq_1]_\ct\\
    &=[p_1,Mq_2]_\ct\\
    &=[Mp_1,q_2]_\ct\\
    &=[Mp_2,q_2]_\ct\\
       &=[Mp_2,Mq_2]_M.
  \end{split}
  \]
  Furthermore, $[\cdot,\cdot]_M$ is not degenerate. Indeed, if for some $p\in\mathbb H_t$ we have $[Mp,Mq]_M=0$ for all $q\in\mathbb H_t$ we obtain
  $[Mp,q]_\ct=0$ for all $q\in\mathbb H_t$ and hence $Mp=0$.
Finally,  $({\rm ran}\, M,[\cdot,\cdot]_M)$ is  a Hilbert space since we are in the finite dimensional setting.\\

STEP 2: {\sl Let $i$ denote the injection from ${\rm ran}\, M$ into $\mathbb H_t^n$. Then $i^{[*]}q=Mq$, where $[*]$ denotes the adjoint between the Pontryagin space
  $\mathbb H_t^n$   and the Hilbert space ${\rm ran}\, M$.}\\

Indeed, for $p,q\in\mathbb H_t^n$ we have
\[
  \begin{split}
    [i(Mp),q]_\ct&=[Mp,q]_\ct\\
    &=[Mp,Mq]_M
  \end{split}
\]
and
\[
  \begin{split}
[i(Mp),q]_\ct&=[Mp,i^{[*]}q]_M.
  \end{split}
\]
Comparing both equalities we obtain the result.\\

STEP 3: {\sl We have $M=ii^{[*]}$.}\\

This follows directly from the previous steps.
  \end{proof}

We can rewrite the above as 
\[
M=iJ_nM^\ct,
\]
where $J_n$ is the extension to $\mathbb H_t^n$ of the symmetry $J$ defined by \eqref{J-pont}.\\

The following result concerns the case $t<0$: 
\begin{theorem}
Assume $t<0$.
  $M\in\mathbb H_t^{n\times n}$ is $\ct$-non-negative if and only if it can be written as $M=FF^\ct$ where $F\in\mathbb  H_t^{n\times n}$, where $F$ may be assumed
  lower triangular.
  \end{theorem}

  \begin{proof}
    The result is true for $n=1$. Then $M=a\ge 0$ and one can take $F=\sqrt{a}$. We proceed by induction and assume $n>1$ and
    the result proved for $n-1$. Let now $M=(m_{jk})_{j,k=1}^n\in \mathbb H_t^{n\times n}$ and assume  $M=M^\ct$.\\

    STEP 1: {\sl Assume $m_{11}=0$. Then, $m_{1k}=0$, $k=2,3,\ldots n$.}\\

   We prove the result for $m_{12}$. The proof for the other indexes is the same. With
\[
    p=\begin{pmatrix}p_1\\ p_2\\ 0\\ \vdots\\ 0\end{pmatrix}
  \]
where $p_1,p_2\in    \mathbb H_t$,  we can write
    \[
   [p,Mp]_\ct={\rm Tr}\begin{pmatrix}p_1^\ct&p_1^\ct\end{pmatrix}\begin{pmatrix}0&m_{12}\\ m_{12}^\ct&m_{22}\end{pmatrix}\begin{pmatrix}p_1\\ p_2\end{pmatrix}
   \ge 0,   \]
 and thus
 \[
{\rm Tr}\, (p_1^\ct m_{12}p_2+p_2^\ct m_{12}^\ct p_1+p_2^\ct m_{22}p_2)\ge 0, \quad \forall p_1,p_2\in    \mathbb H_t.
   \]
   Write $p_1=a+j_tb$ and $m_{12}=c+j_td$. Moreover let $p_2=\e>0$. Using \eqref{sym-345} we obtain
   \[
a\overline{c}+\overline{a}c-t(b\overline{d}+\overline{b}d)\e\ge0.
\]
 Letting $\e$ go to $0$ we obtain
\[
a\overline{c}+\overline{a}c\ge0,\quad \forall a\in\mathbb C.
  \]
  The choice $a=1$ and $a=-1$ lead to ${\rm Re}\,c=0$ and the choices $a=i$ and $a=-i$ gives ${\rm Im}\, c=0$. We then have
  \[
-t(b\overline{d}+\overline{b}d)\e\ge0,\quad \forall b\in\mathbb C.
\]
Setting now $b\in\mathbb R$ we get
\[
-2tb(d+\overline{d})\ge 0,\quad \forall b\in\mathbb R
\]
which cannot be unless $d+\overline{d}=0$ since $t\not=0$. A similar argument shows that ${\rm Im}\,d=0$ and hence $d=0$.\\
     
STEP 2: {\sl  Reiterating Step 1 a number of times, we get  $M=0_{\mathbb H_t^{N\times N}}$ unless there exists $n_0\in\left\{1,\ldots, n\right\}$
  such that $m_{n_0n_0}>0$.}\\

Indeed, if $m_{22}>0$ the step is proved. If $m_{22}=0$, the preceding argument shows that $m_{22}=m_{23}=\cdots=m_{2N}=0$, and we already have that $m_{12}=0$.
Iterating this argument leads to the required result.\\

    STEP 3:  {\sl Assume that $m_{n_0n_0}\not=0$ but $m_{kk}=0$, $k=1,\ldots n_0-1$ (if $n_0=1$ we mean $m_{11}\not=0$). Then
    \[
      M=\begin{pmatrix}0_{\mathbb H_t^{(n_0-1)\times (n_0-1)}}&0_{\mathbb H_t^{(n_0-1)\times (n-n_0+1)}}\\
        0_{\mathbb H_t^{(n-n_0+1)\times (n_0-1)}}&T\end{pmatrix}
      \]
      where
      \[
        T=\begin{pmatrix} x_{11}&B\\
          B^{\ct}&D\end{pmatrix}
        \]
        with $x_{11}>0$.}\\

      This is a consequence of Step 1 and Step 2.\\

      STEP 4: {\sl $T$ is $\ct$-non-negative and $x_{11}>0$.}\\

      $x_{11}>0$ by construction. Furthermore, with
      \[
        p=\begin{pmatrix}0_{\mathbb H_t^{(n_0-1)}}\\ q\end{pmatrix}, \quad q\in\mathbb H_t^{(n-n_0+1)},
      \]
      we have
      \[
        [p,Mp]_\ct=[q,Tq]_{\ct}\ge 0,\quad \forall q\in\mathbb H_t^{(n-n_0+1)},
      \]
      and so $T\ge_{\ct} 0$.\\
      
STEP 5: {\sl $T$ can be factorized as
        \[
          T=\begin{pmatrix}1&0\\B^{\ct}x_{11}^{-1}&I\end{pmatrix}\begin{pmatrix}x_{11}&0\\ 0
            &D-B^\ct x_{11}^{-1}B\end{pmatrix}\begin{pmatrix}1&x_{11}^{-1}B\\0&I\end{pmatrix}.
        \]
      }\\

      This is a direct computation, based on a well-known factorization formula for block matrices.\\
      
      STEP 6: {\sl In the notation of Step 3, $D-B^\ct x_{11}^{-1}B\ge_\ct 0$.}\\

      Define
      \[
     q=\begin{pmatrix}1&x_{11}^{-1}B\\0&I\end{pmatrix}^{-1}r=\begin{pmatrix}1&-x_{11}^{-1}B\\0&I\end{pmatrix}r,
   \]
   where $r\in\mathbb H_t^{(n-n_0)}$. Then,
   \[
     [r, (D-B^\ct x_{11}^{-1}B)r]_\ct=[p,Tp]_\ct \ge_\ct0,\quad\forall r\in\mathbb H_t^{(n-n_0)}.
   \]
   
   STEP 7: {\sl Apply the hypothesis induction to $D-B^\ct x_{11}^{-1}B$ to conclude the proof.}\\

   We write
   \[
     D-B^\ct x_{11}^{-1}B=RR^\ct,\quad{\rm where}\quad R\in\mathbb H_t^{(n-n_0)\times(n-n_0)}.
   \]
   We have
   \[
     \begin{split}
     M&=\begin{pmatrix}0_{\mathbb H_t^{(n_0-1)\times (n_0-1)}}&0_{\mathbb H_t^{(n_0-1)\times (n-n_0+1)}}\\
       0_{\mathbb H_t^{(n-n_0+1)\times (n_0-1)}}&T\end{pmatrix}\\
     &=  \begin{pmatrix}0_{\mathbb H_t^{(n_0-1)\times (n_0-1)}}&0_{\mathbb H_t^{(n_0-1)\times (n-n_0+1)}}\\
  0_{\mathbb H_t^{(n-n_0+1)\times (n_0-1)}}&
  \begin{pmatrix}1&0\\B^{\ct}x_{11}^{-1}&I\end{pmatrix}\begin{pmatrix}x_{11}&0\\ 0
    &D-B^\ct x_{11}^{-1}B\end{pmatrix}\begin{pmatrix}1&x_{11}^{-1}B\\0&I\end{pmatrix}\end{pmatrix}\\
& =
\begin{pmatrix}0_{\mathbb H_t^{(n_0-1)\times (n_0-1)}}&0_{\mathbb H_t^{(n_0-1)\times (n-n_0+1)}}\\
  0_{\mathbb H_t^{(n-n_0+1)\times (n_0-1)}}&
  \begin{pmatrix}1&0\\B^{\ct}x_{11}^{-1}&I\end{pmatrix}
  \begin{pmatrix}x_{11}&0\\ 0
            &RR^\ct\end{pmatrix}\begin{pmatrix}1&x_{11}^{-1}B\\0&I\end{pmatrix}\\
        \end{pmatrix}
      \end{split}
      \]
      But
      \[
\begin{split}
  \begin{pmatrix}1&0\\B^{\ct}x_{11}^{-1}&1\end{pmatrix}
  \begin{pmatrix}x_{11}&0\\ 0
          &RR^\ct\end{pmatrix}\begin{pmatrix}1&x_{11}^{-1}B\\0&I\end{pmatrix}=\\
        &\hspace{-4cm}=
        \begin{pmatrix}1&0\\B^{\ct}x_{11}^{-1}&I\end{pmatrix}\begin{pmatrix}\sqrt{x_{11}}&0\\0&R\end{pmatrix}
        \begin{pmatrix}\sqrt{x_{11}}&0\\0&R^\ct\end{pmatrix}\begin{pmatrix}1&x_{11}^{-1}B\\0&I\end{pmatrix}\\
        &\hspace{-4cm}=\begin{pmatrix}\sqrt{x_{11}}&0\\ B^\ct x_{11}^{-1/2}&R\end{pmatrix}\begin{pmatrix}\sqrt{x_{11}}
          &B x_{11}^{-1/2}\\ 0&R^\ct\end{pmatrix}.
              \end{split}
    \]
    Thus
\[
M=        
\underbrace{\begin{pmatrix}0_{\mathbb H_t^{(n_0-1)\times (n_0-1)}}&0_{\mathbb H_t^{(n_0-1)\times (n-n_0+1)}}\\
  0_{\mathbb H_t^{(n-n_0+1)\times (n_0-1)}}&
  \begin{pmatrix}1\sqrt{x_{11}}&0\\B^{\ct}x_{11}^{-1/2}&R\end{pmatrix}
\end{pmatrix}}_{F}
\underbrace{\begin{pmatrix}0_{\mathbb H_t^{(n_0-1)\times (n_0-1)}}&0_{\mathbb H_t^{(n_0-1)\times (n-n_0+1)}}\\
  0_{\mathbb H_t^{(n-n_0+1)\times (n_0-1)}}&
  \begin{pmatrix}1\sqrt{x_{11}}&B x_{11}^{-1/2}\\ 0&R^\ct\end{pmatrix}
\end{pmatrix}}_{F^\ct},
  \]
  and so $M$ has the required factorization.\\

  STEP 8: {\sl Any matrix of the form $FF^\ct$ is $\ct$-non-negative.}\\

  It suffices to write
  \[
    [p,FF^\ct p]_\ct=[F^\ct p,F^\ct p]_\ct\ge 0.
    \]
\end{proof}


\begin{proposition}
Assume $t<0$ and $A=A^\ct$. The eigenvalues of $A$ are then real.
  \end{proposition}

  \begin{proof}
    We define on $\mathbb H_t^n$ the $\mathbb H_t$-valued form
    \begin{equation}
\left\{f,g\right\}_\ct=g^{\ct}f.
\end{equation}
In particular we have
\begin{equation}
  \label{matrix-valued}
  \left\{f,f\right\}=(|a|^2-t|b|^2) I_2.
  \end{equation}
Note that, for $f,g\in\mathbb H_t^n$, $c\in\mathbb H_t$  and $A\in\mathbb H_t^{n\times n}$,
\[
  \left\{f,gc\right\}_\ct=c^\ct g^{\ct}f=c^\ct\left\{f,g\right\}\quad{\rm and}\quad
  \left\{fc,g\right\}_\ct g^{\ct}fc.=c^\ct\left\{f,g\right\}c.
  \]
Assume now that $Af=f\lambda$ where $\lambda\in\mathbb C$. In view of the above, we have
\[
  \left\{Af,f\right\}=\begin{cases}\,\,=\, \left\{f\lambda,f\right\}\,\, =\,\, \left\{f,f\right\}\lambda\\
        \,\,=\, \left\{f,Af\right\}\,\, =\,\, \left\{f,fc\right\}\,\, =\,\,\overline{\lambda}\left\{f,f\right\}.
\end{cases}
\]
since $\lambda^\ct=\overline{\lambda}$ for $\lambda\in\mathbb C$. It follows that
\[
  \left\{f,f\right\}\lambda=\overline{\lambda}\left\{f,f\right\}.
  \]
Equation \eqref{matrix-valued} allows to conclude since $t<0$.
\end{proof}

\begin{corollary}
Assume that $A=A^\ct$. Then eigenvectors corresponding to different eigenvalues are orthogonal in the $[\cdot,\cdot]_\ct$ form.
\end{corollary}

\begin{proof}
Let $f,g\in\mathbb H_t^n$, corresponding to the eigenvalues $\lambda$ and $\mu$ respectively. Then, since $\lambda $ and $\mu$ are real, 
\[
  [Af,g]_\ct=[f\lambda,g]_\ct=\lambda [f,g]_\ct
\]
on the one hand, and
\[
  [Af,g]_\ct=[f,Ag]_\ct=\mu [f,g]_\ct
\]
on the other hand. It follows that
\[
(\lambda-\mu)[f,g]_\ct=0
\]
and hence the result.
  \end{proof}
\subsection{$\ct$-positive definite kernels}
The $\ct$-adjoint does not behave well with respect to the classical adjoint (see Remark \ref{remstar}), and the fact that we have real Hilbert spaces here rather than
Hilbert spaces over the complex numbers play an important role, both in the definitions and in the arguments.

\begin{definition}
Let $\Omega$ be some set.    
The $\mathbb H_t^{n\times n}$-valued function $K(z,w)$ defined for $z,w\in \Omega$ is called $\ge_\ct$-positive definite if it is $\ct$-Hermitian:
\begin{equation}
  K(z,w)=K(w,z)^\ct,\quad\forall w,z\in \Omega,
\end{equation}
and if for every $N\in\mathbb N$, any $w_1,\ldots, w_N\in \Omega$ and $c_1,\ldots, c_N\in\mathbb H_t^N$,
\begin{equation}
  \label{pos1}
{\rm Tr}\,\left(\sum_{j,k=1}^Nc_k^\ct K(w_k,w_j)c_j\right)\ge 0.
  \end{equation}
\end{definition}

\begin{theorem}
  Let $K(z,w)$ be a $\mathbb H_t^{n\times n}$-valued function $\ct$-positive definite in $\Omega$. Then there exists a uniquely defined real Hilbert right $\ct$-module
  $\mathfrak H(K)$ of functions defined on $\Omega$ with
  reproducing kernel $K(z,w)$, meaning that for every $w\in\Omega$ and $c\in\mathbb H_t^n$  the function
  \[
z\mapsto K(z,w)c\in\mathfrak H(K),
\]
and that for every $f\in\mathfrak H(K)$,
\begin{equation}
  \label{rk}
  [f,K(\cdot,w)c]_K={\rm Tr}\, [f(w),c]_{\ct},
\end{equation}
where $[\cdot,\cdot]_K$ is the real-valued inner product on $\mathfrak H(K)$.
\end{theorem}

{\bf Proof:} The proof follows closely the classical complex setting case. We repeat the arguments, suitably adapted to the present framework, and proceed in a number
of steps. We set
$\stackrel{\circ}{\mathfrak H}\hspace{-1mm}(K)$ to be the linear space of the functions of the form
\[
  f(z)=\sum_{k=1}^m K(z,w_k)c_k,\quad m=1,2,\ldots,\quad w_1,\ldots, w_m\in \Omega,\quad c_1,\ldots, c_m\in\mathbb H_t^n
  \]
  with real-valued symmetric form
  \begin{equation}
    \label{inner-prod}
[f,g]_K={\rm Tr}\,\left(\sum_{k=1}^m\sum_{j=1}^r [K(v_j,w_k)c_k,d_j]_\ct\right),
\end{equation}
where $g(z)=\sum_{j=1}^r K(z,v_j)d_j,\quad r\in \mathbb N,\quad v_1,\ldots, v_m\in \Omega,\quad d_1,\ldots, d_r\in\mathbb H_t^n$. The representations of $f$ and $g$
need not be unique, and we begin  with:\\

STEP 1: {\sl \eqref{inner-prod} is well defined, i.e. does not depend on the given representations of $f$ and $g$ as linear combination of kernels.}\smallskip

The proof is the same as for the classical setting. See e.g. \cite{aron}.\\

STEP 2: {\sl $[\cdot,\cdot]_K$ is a positive symmetric non-degenerate form and the reproducing property \eqref{rk} holds in $\stackrel{\circ}{\mathfrak H}\hspace{-1mm}(K)$.}\smallskip

Indeed, setting $g(\cdot)=K(\cdot, w)c$ in \eqref{inner-prod} we obtain \eqref{rk}. This equality implies in turn that $[\cdot,\cdot]_K$ is non-degenerate. If $[f,g]_K=0$ for all
$g\in\stackrel{\circ}{\mathfrak H}\hspace{-1mm}(K)$, the choice $g(\cdot)=K(\cdot, w)c$ shows that
\[
  {\rm Tr}\,[f(w),c]_\ct=0,\quad \forall c\in\mathbb H_t^N,\quad \forall w\in\Omega,
  \]
  and so (since $[\cdot,\cdot]_\ct$ is non-degenerate)
\[
c^\ct f(w)=0,\quad \forall c\in\mathbb H_t^N,\quad \forall w\in\Omega,
\]
and hence $f=0$. The positivity of the form follows from \eqref{pos1}. We therefore also have proved Step 3 below.\\

STEP 3: {\sl The pair $(\stackrel{\circ}{\mathfrak H}\hspace{-1mm}(K),[\cdot,\cdot]_K)$ is a real pre-Hilbert space}\smallskip

If the real linear span of the functions $K(\cdot, w)c$ with $w\in\Omega$ and $c\in\mathbb H_t^n$, is finite dimensional the proof is finished. The following step deals with the
case when the dimension is not finite. By a standard result in metric space theory, there is a uniquely defined, up to a metric space isomorphism, Hilbert space in which
$\stackrel{\circ}{\mathfrak H}\hspace{-1mm}(K)$ is dense.

STEP 4: {\sl  One can choose the closure of $\stackrel{\circ}{\mathfrak H}\hspace{-1mm}(K)$ to be a space of functions on $\Omega$, which we denote by
  $\mathfrak H(K)$, and the reproducing kernel property holds in $\mathfrak H(K)$.}\smallskip

To construct the completion one looks at the family of Cauchy sequences $(f_n)$ in $\stackrel{\circ}{\mathfrak H}\hspace{-1mm}(K)$,
and define an equivalence class as follows: two Cauchy sequences $(f_n)$ and $g_n)$ will be said to be equivalent if
\[
\lim_{n\rightarrow\infty} [f_n-g_n,f_n-g_n]_K=0.
\]
This indeed defines an equivalence relation, and we denote the equivalence classes as $f=\widetilde{(f_n)}$.
The space of such equivalence classes with the inner product
\begin{equation}
  \label{wsx}
  [f,g]=\lim_{n\rightarrow\infty}[f_n,g_n]_K,\quad (f_n)\in f,\quad (g_n)\in g
  \end{equation}
  (which is well defined) defines a real Hilbert space. The identification with a space of functions on $\Omega$ is then done as follows:
  To such an equivalence class one associates the function
\[
  f(w)=\lim_{n\rightarrow\infty} f_n(w),
\]
and the limit does not depend on the chosen representative $(f_n)$ of the equivalence class.
The space ${\mathfrak H}(K)$ consists exactly of these limit functions, with the induced inner product \eqref{wsx}
and has the reproducing kernel property.
\mbox{}\qed\mbox{}\\

  \section{Realization theory in $\mathbb H_t$}
  \setcounter{equation}{0}
  \label{sec:real}
\subsection{Definitions}
In our previous paper \cite{acv1} we introduced in the present setting the notions of rational functions, one of the equivalent definitions being in terms of a realization, namely:
\begin{definition}
  The $\mathbb H_t^{m\times n}$-valued function $R(x)$ of the real variable $x$ defined for $x=0$ is rational if it can be written as
  \begin{equation}
    R(x)=D+xC(I_{\mathbb H_t^N}-xA)^{-1}B
    \end{equation}
    with $D=R(0)$ and matrices $(C,A,B)\in\mathbb H_t^{m\times N}\times\mathbb H_t^{N\times N}\times\mathbb H_t^{N\times n}$.
  \end{definition}

We here recall some formulas, well known in the complex setting. The proofs are the same and omitted; see for instance \cite{bgk1,ot_1_new}.

\begin{proposition}
  The product of two rational functions of compatible sizes with realizations $R_j(x)=D_j+xC_j(I_{\mathbb H_t^{N_j}}-xA_j)^{-1}B_j)$, $j=1,2$,
  has a realization, given by
  \begin{equation}
  \label{prod_1}
  \begin{split}
    (D_1+xC_1(I_{\mathbb H_t^{N_1}}-xA_1)^{-1}B_1)(D_2+xC_2(I_{\mathbb H_t^{N_2}}-xA_2)^{-1}B_2)&=\\
    &\hspace{-20mm}=D+xC(I_{\mathbb H_t^{N}}-xA)^{-1}B,
    \end{split}
\end{equation}
with
\begin{eqnarray}
  \label{prod_2}
  N&=&N_1+N_2\\
  \label{prod_3}
  A&=&\begin{pmatrix}A_1&B_1C_2\\
    0_{\mathbb H_t^{N_1}}&A_2\end{pmatrix}\\
  B&=&\begin{pmatrix}B_1D_2\\ B_2\end{pmatrix}\\
  \label{prod_0}
  C&=&\begin{pmatrix}C_1& D_1C_2\end{pmatrix}\\
  D&=&D_1D_2
       \label{prod_6}
\end{eqnarray}
\end{proposition}
\begin{proposition}
  Assume $R(x)$ rational $\mathbb H_t^{n\times n}$-valued and $R(0)$-invertible, with realization
  \[
R(x)=D+xC(I_{\mathbb H_t^N} -xA)^{-1}B.
    \]
    Then,
    \begin{equation}
      \label{inverse11}
      R(x)^{-1}=D^{-1}-xD^{-1}C(I_{\mathbb H_t^N} -xA^\times)^{-1}BD^{-1}
    \end{equation}
    with
    \begin{equation}
    \label{a-times}
A^\times=A-BD^{-1}C,
\end{equation}
  \end{proposition}

  We note that $A^\times$ is the Schur complement of $A$ in the block matrix representation $\begin{pmatrix}A&B\\ C&D\end{pmatrix}$ of the realization.
In the realization \eqref{prod_3}-\eqref{prod_6}, we also note that
\begin{equation}
  \label{atimes}
  A^\times =A-BD^{-1}C=\begin{pmatrix}A_1-B_1D_1^{-1}C_1&0\\    -B_2D_2^{-1}D_1^{-1}C_1&A_2-B_2D_2^{-1}C_2\end{pmatrix}
\end{equation}
and is block lower triangular while $A$ is block  upper triangular.
This is a key computation in the description of the minimal factorizations; see Section \ref{mini-facto}.
      \subsection{Observability}
  \begin{definition}
    The pair $(G,T)\in\mathbb C^{a\times M}\times\mathbb C^{M\times M}$ is observable if the following condition is in force:
    Let $\xi\in\mathbb C^M$. Then,
    \[
   G(I_{M}-xT)^{-1}\xi\equiv0,\quad x\in(-\e,\e) \quad \mbox{for some $\e>0$}\quad \iff \xi=0.
  \]
  The pair $(C,A)\in\mathbb H_t^{n\times N}\times\mathbb H_t^{N\times N}$ is observable if the following condition is in force:
      Let $h\in\mathbb H_t^N$. Then,
    \[
   C(I_{\mathbb H_t^N}-xA)^{-1}h\equiv0,\quad x\in(-\e,\e) \quad \mbox{for some $\e>0$}\quad \iff h=0.
  \]
  \end{definition}
  
  \begin{proposition}
    \label{C-A-I}
The pair $(C,A)\in\mathbb H_t^{n\times N}\times\mathbb H_t^{N\times N}$ is observable if and only if the pair $(I(C),I(A))\in\mathbb C^{2n\times 2N}\times\mathbb C^{2N\times 2N}$ is observable.
\end{proposition}

\begin{proof}
  Assume that $(C,A)$ is observable and let $\xi\in\mathbb C^{2N}$ be such that for some $\e>0$
  \[
    I(C)(I_{2N}-xI(A))^{-1}\xi\equiv0,\quad x\in(-\e,\e).
  \]
  Write
  \[
    \xi=\begin{pmatrix}a_1\\ \vspace{-3mm}\\
      \overline{b_1}\\a_2\\ \vspace{-3mm}\\\overline{b_2}\\ \vdots \\a_N\\ \vspace{-3mm}\\\overline{b_N}\end{pmatrix}
\]
and define $\eta\in\mathbb C^{2N\times 2}$ by
\[
  \eta=\begin{pmatrix}a_1&tb_1 \\ \vspace{-3mm}\\ \overline{b_1}&\overline{a_1}\\a_2& tb_2\\ \vspace{-3mm}\\ \overline{b_2}&\overline{a_2}\\ \vdots&\vdots\\ a_N&tb_N\\ \vspace{-3mm}\\ \overline{b_N}&\overline{a_N}\end{pmatrix}.
\]
Then, by Lemma \ref{debut0},
  \[
    I(C)(I_{2N}-xI(A))^{-1}\eta\equiv0,\quad x\in(-\e,\e).
  \]
  and so
    \[
    I(C)(I_{2N}-xI(A))^{-1}\begin{pmatrix}\xi&\eta\end{pmatrix}\equiv0,\quad x\in(-\e,\e).
  \]
  But
\[
  \eta=I(h)\quad {\rm where}\quad h=\begin{pmatrix}  a_1+b_1j_t\\ a_2+b_2 j_t\\ \vdots \\a_N+b_Nj_t \end{pmatrix}.
\]
and so we have by Lemma \ref{lemm-ab-I}
\begin{equation}
  \label{c-a-obs}
C(I_{\mathbb H_t^N}-xA)^{-1}h\equiv0,\quad x\in(-\e,\e).
  \end{equation}
  so that $h=0$ since $(C,A)$ is observable. Hence $\xi=0$ and so $(I(C),I(A))$ is observable.\smallskip

  Assume now that $(I(C),I(A))$ is observable and let $h\in\mathbb H_t^N$ be such that \eqref{c-a-obs} is in force. Then,
    \[
      I(C)(I_{2N}-xI(A))^{-1}I(h)=0,\quad x\in(\e,\e),
    \]
    and since the pair $(I(C),I(A))$ is observable we have $\xi=\eta=0$ where $I(h)=\eta$, and so $h=0$.
  \end{proof}
  \subsection{Controllability}
  \begin{definition}
    The pair $(T,F)\in\mathbb C^{M\times M}\times\mathbb C^{M\times b}$ is controllable if the following condition is in force:
    \begin{equation}
      \cup_{k=0}^\infty{\rm ran}\, T^kF=\mathbb C^M.
    \end{equation}
  \end{definition}

Note that by the Cayley-Hamilton theorem we have
\[
  \cup_{k=0}^\infty{\rm ran}\, T^kF=\cup_{k=0}^{N-1}{\rm ran}\, T^kF
\]
for some $N\le M\in\mathbb N$.
    \begin{lemma}
      Let $(T,F)\in\mathbb C^{M\times M}\times\mathbb C^{M\times b}$. The following are equivalent:\\
      $(1)$ The pair $(T,F)$ is controllable.\\
      $(2)$ The pair $(F^*,T^*)$ is observable.\\
      $(3)$ If $\xi\in\mathbb C^{1\times M}$ is such that
      \[
        \xi(I_M-xT)^{-1}F=0_{1\times b}
        \]
      in a real neighborhood of $0$, then $\xi=0_{1\times M}$.\\
      $(4)$
      \begin{equation}
        \label{non-dense}
\xi T^kF=0_{1\times b},\quad k=0,1,2,\ldots\quad\iff\quad \xi=0_{1\times M}.
        \end{equation}
      \end{lemma}

See for instance \cite{bgk1,ot_1_new} for a proof.\\

      We remark that $\xi^*\in\mathbb C^M$ and that \eqref{non-dense} is equivalent to
      \begin{equation}
        \langle T^kF\zeta,\xi^*\rangle=0,\quad k=0,1,\ldots,\,\,\forall\zeta\in\mathbb C^m,
      \end{equation}
      where the brackets denote the standard inner product in $\mathbb C^m$.\smallskip

      \begin{definition}
        The pair $(A,B)\in\mathbb H_t^{N\times N}\times\mathbb H_t^{N\times m}$ is controllable if the following condition is in force:
        Let $h\in\mathbb H_t^{1\times N}$. Then
        \begin{equation}
          \label{contro-4}
          h A^kB=0_{1\times m},\quad k=0,1,\ldots\quad \iff\quad h=0_{\mathbb H_t^{1\times N}},
        \end{equation}
        or equivalently
        \begin{equation}
          h(I_{\mathbb H_t^N}-xA)^{-1}B=0_{1\times m},\quad x\in(-\e,\e)\quad \iff h=0_{\mathbb H_t^{1\times N}}.
          \end{equation}
        \end{definition}

      The equivalence of the two conditions in this definition follows from the series expansion at the origin
      \[
(I_{\mathbb H_t^N}-xA)^{-1}=\sum_{\ell=0}^\infty\e^\ell A^n,
\]
where the convergence can be understood in $\mathbb C^{2N\times 2N}$ by using the map $I$.
        
\begin{proposition}
  \label{B-A-I}
      The pair $(A,B)\in\mathbb H_t^{N\times N}\times\mathbb H_t^{N\times m}$ is controllable if and only if the pair $(I(A),I(B))\in\mathbb C^{2N\times 2N}\times
      \mathbb C^{2N\times 2m}$ is controllable.
\end{proposition}

\begin{proof}
  Assume that $(A,B)$ is controllable, and let $\xi\in\mathbb C^{1\times 2N}$ be such that for some $\e>0$
  \[
\xi  (I_{2N}-xI(A))^{-1}I(B)\equiv0,\quad x\in(-\e,\e).
\]
Write
\[
   \xi=\begin{pmatrix}a_1&   tb_1&a_2&tb_2& \cdots & a_N&tb_N\end{pmatrix}
\]
and define $\eta\in\mathbb C^{2\times 2N}$ by
\[
  \eta=\begin{pmatrix}\overline{b_1}&\overline{a_1}&\overline{b_2}&\overline{a_2}&\cdots&\overline{b_N}&\overline{a_N}\end{pmatrix}.
\]
Then, by Lemma \ref{debut0},
  \[
   \eta(I_{2N}-xI(A))^{-1}I(B)\equiv0,\quad x\in(-\e,\e),
 \]
 and so
   \[
\begin{pmatrix}\xi\\\eta\end{pmatrix}     (I_{2N}-xI(A))^{-1}I(B)     \equiv0,\quad x\in(-\e,\e).
  \]
  But
  \[
 \eta=I(h)\quad {\rm where}\quad h=\begin{pmatrix}  a_1+b_1j_t& a_2+b_2 j_t& \cdots &a_N+b_Nj_t \end{pmatrix}.\]
and so we have by Lemma \ref{lemm-ab-I}
\begin{equation}
\label{c-a-contro}
h(I_{\mathbb H_t^N}-xA)^{-1}B\equiv0,\quad x\in(-\e,\e),
 \end{equation}
 so that $h=0$ since $(A,B)$ is controllable. Hence $\xi=0$ and so $(I(A),I(B))$ is controllable.\smallskip

Assume now that $(I(A),I(B))$ is controllable and let $h\in\mathbb H_t^{1\times N}$ be such that \eqref{c-a-contro} is in force. Then,
   \[
     I(h)(I_{2N}-xI(A))^{-1}I(B)=0,\quad x\in(\e,\e),
    \]
    and since the pair $(I(A),I(B))$ is controllable we have $\xi=\eta=0$ where $I(h)=\begin{pmatrix}\xi
&\eta\end{pmatrix}$.
\end{proof}

\begin{proposition}
  The following are equivalent:\\
$(1)$ The pair $(A,B)\in\mathbb H_t^{N\times N}\times\mathbb H_t^{N\times n}$ is controllable.\\
$(2)$ The pair $(B^{\circledast},A^{\circledast})\in\mathbb H_t^{m\times N}\times\mathbb H_t^{N\times N}$ is observable.\\
\end{proposition}

\begin{proof}
  Assume that $(1)$ is in force and let $h\in\mathbb H_t^N$ be such that
  \[
    B^{\circledast}(I_{\mathbb H_t^N}-xA^{\circledast})^{-1}h=0,\quad x\in(-\e,\e).
  \]
  Since $\circledast$ is contravariant we have
  \[
    h^{\circledast}(I_{\mathbb H_t^N}-xA)^{-1}B=0,\quad x\in(-\e,\e).
  \]
  Since the pair $(A,B)$ is controllable, we have $h^{\circledast}=0$ and since $\circledast$ is involutive, $h=0$; thus $(2)$ holds. The proof that $(2)$ implies $(1)$ goes in
  the same way.
\end{proof}

\begin{lemma}
  The pair $(A,B)\in\mathbb H_t^{N\times N}\times\mathbb H_t^{N\times m}$ is controllable if and only if the linear span $\mathcal M$ of the vectors $(I_{\mathbb H_t^{N}}-yA)^{-1}Bb$ with $y\in(\e,\e)$ and $b$ runs
  through $\mathbb H_t^m$ is all of $\mathbb H_t^N$.
  \end{lemma}
  \begin{proof} 
    $\mathbb H_t^N$ with the form \eqref{inner-prod-h-t} is a finite dimensional Pontryagin space with indefinite inner product $[\cdot,\cdot]_{\circledast}$. Let $J$ be a signature such that $[\cdot, J\cdot]_\ct$ is a Hilbert space. Then the above span will not be equal to $\mathbb H_t^N$ if and only if there is
    a $h\in\mathbb H_t^N$ such that
    $[h,Jv]=0$ for $v$ in $\mathcal M$. So $[Jh,v]=0$ for all $v\in\mathcal M$.
  \end{proof}

  \begin{lemma}
    \eqref{contro-4} is equivalent to
    \begin{equation*}
          [h^{\circledast}, A^kBb]_{\circledast}=0,\quad k=0,1,\ldots,\quad b\in\mathbb H_t^m \iff\quad h=0_{\mathbb H_t^{1\times N}}.
\end{equation*}
        \end{lemma}

        \begin{proof}
          Assume that ${\rm span}\left\{A^kBb,\quad k=0,1,\ldots,\,\,b\in\mathbb H_t^m \right\}=\mathbb H_t^N$. Since
          $(\mathbb H_t^N,[\cdot,\cdot]_{\circledast})$
          is a real Pontryagin space (see Proposition \ref{pontry}), and in particular non-degenerate. Thus $h^{\circledast}=0_{\mathbb H_t^N}$ and in particular $h A^kB=0$ for $k=0,1,\ldots$.\smallskip

          Conversely, if $h A^kB=0$ for $k=0,1,\ldots$, then ${\rm Re}\, h A^kBb=0$ for $k=0,1,\ldots$. But ${\rm Re}\, h A^kBb=[h^{\circledast},
            A^kBb]_{\circledast}$, and $h=0$ by the previous arguments.
\end{proof}

\subsection{Minimal realizations}
\label{minimal-real}
We note that a realization is never unique since one can write
\[
r(x)=D+xC(I_{\mathbb H_t^N}-xA)^{-1}B=D+xCS^{-1}(I_{\mathbb H_t^N}-xS^{-1}AS)^{-1}SB
  \]
where $S\in\mathbb H_t^{N\times N}$ is invertible in $\mathbb H_t^{N\times N}$. When the realization is minimal this is the only degree of freedom.

\begin{theorem}
  Let $r(x)=D+C_1(I_{\mathbb H_t^{N_1}}-xA_1)^{-1}B_1=D+C_2(I_{\mathbb H_t^{N_2}}-xA_2)^{-1}B_2$ be two minimal realizations of the $\mathbb H_t^{n\times m}$-valued rational function $r(x)$ where in both cases $D=r(0)$. Then, $N_1=N_2$ and the two realizations are similar, meaning that there exists
  an invertible matrix $S\in\mathbb H_t^{N\times N}$ such that
  \begin{equation}
    \begin{pmatrix}S&0\\ 0&I_{\mathbb H_t^n}\end{pmatrix}
\begin{pmatrix}A_1&B_1\\C_1&D\end{pmatrix}=\begin{pmatrix}A_2&B_2\\C_2&D\end{pmatrix}.
        \begin{pmatrix}S&0\\ 0&I_{\mathbb H_t^n}\end{pmatrix}
\end{equation}
  Furthermore, $S$ is uniquely defined by the two realizations.
  \end{theorem}

  \begin{proof}
    We follow the proof described in \cite[Theorem 9.2.4 p. 247]{zbMATH06658818} for the quaternionic setting. A key fact which allows to adapt the proof
    is that, since $\mathbb H_t^{N_1}$ and $\mathbb H_t^{N_2}$ are $\mathbb H_t$ modules with basis, linear maps are represented by matrices.\smallskip

As in \cite[(9.15) p. 247]{zbMATH06658818} we first write that
    \begin{eqnarray}
      \label{cab-1}
      \frac{r(x)-r(y)}{x-y}&=&C_1(I_{\mathbb H_t^{N_1}}-xA_1)^{-1}(I_{\mathbb H_t^{N_1}}-yA_1)^{-1}B_1\\
                           &=&C_2(I_{\mathbb H_t^{N_2}}-xA_2)^{-1}(I_{\mathbb H_t^{N_2}}-yA_2)^{-1}B_2,\quad x,y\in(-\e,\e).
                               \label{cab-2}
      \end{eqnarray}

      We then proceed in a number of steps.\\

      STEP 1: {\sl The ({\sl a priori} multi-valued) formulas}
      \begin{eqnarray}
        \label{U-formula}
          U\left(\sum_{e=1}^E(I_{\mathbb H_t^{N_1}}-y_eA_1)^{-1}B_1a_e\right)
          &=&\sum_{e=1}^E(I_{\mathbb H_t^{N_2}}-y_eA_2)^{-1}B_2a_e\\
          V\left(\sum_{e=1}^E(I_{\mathbb H_t^{N_1}}-x_eA_1^{\circledast})^{-1}C_1^{\circledast}a_e\right)&=
            &\sum_{e=1}^E(I_{\mathbb H_t^{N_2}}-x_eA_2^{\circledast})^{-1}C_2^{\circledast}a_e
              \label{v-formula}
        \end{eqnarray}
        {\sl induce uniquely defined maps $U$ and $V$ from $\mathbb H_t^{N_1}$ into $\mathbb H_t^{N_2}$.}\smallskip

        We first show that in both cases, the function $0_{\mathbb H_t^{N_1}}$ goes to the function $0_{\mathbb H_t^{N_1}}$, and so $U$ and $V$ are
        linear maps. 
        We then show that their domain is $\mathbb H_t^{N_1}$. For $U$, assume
        \[
\sum_{e=1}^E(I_{\mathbb H_t^{N_1}}-y_eA_1)^{-1}B_1a_e=0_{\mathbb H_t^{N_1}}.
\]
Multiplying on the left by $C_1(I_{\mathbb H_t^{N_1}}-xA_1)^{-1}$ and using \eqref{cab-1} and \eqref{cab-2} we get
\[
  \begin{split}
    C_1(I_{\mathbb H_t^{N_1}}-xA_1)^{-1}\left(\sum_{e=1}^E(I_{\mathbb H_t^{N_1}}-y_eA_1)^{-1}B_1a_e\right)&=\sum_{e=1}^E\frac{r(x)-r(y_e)}{x-y_e}a_e\\
    &\hspace{-12mm}=C_2(I_{\mathbb H_t^{N_2}}-xA_2)^{-1}\left(\sum_{e=1}^E(I_{\mathbb H_t^{N_2}}-y_eA_2)^{-1}B_2a_e\right)\\
    &\hspace{-12mm}=0_{\mathbb H_t^n}.
   \end{split}
 \]
 Since the pair $(C_2,A_2)$ is observable it follows that $\sum_{e=1}^E(I_{\mathbb H_t^{N_2}}-y_eA_2)^{-1}B_2a_e=0$.  That the domain of $U$ is all
 of $\mathbb H_t^{N_1}$ follows from the controllability of the pair $(A_1,B_1)$.\smallskip

 The argument for $V$ is similar and goes as follows: Assume
 \[
   \sum_{e=1}^E(I_{\mathbb H_t^{N_1}}-x_eA_1^{\circledast})^{-1}C_1^{\circledast}a_e=0_{\mathbb H_t^{N_1}}
 \]
 Multiplying on the left by $B_1^{\circledast}(I_{\mathbb H_t^{N_1}}-yA_1^{\circledast})^{-1}$  and applying $\circledast$ on both sides we get
 \[
\left(\sum_{e=1}^Ea_e^\circledast C_1(I_{\mathbb H_t^{N_1}}-x_eA_1)^{-1}\right)(I_{\mathbb H_t^{N_1}}-yA_1)^{-1}B_1=0,\quad y\in(-\e,\e).
   \]

 and using \eqref{cab-1} and \eqref{cab-2} we obtain
 \[
  \left(\sum_{e=1}^Ea_e^\circledast C_2(I_{\mathbb H_t^{N_2}}-x_eA_2)^{-1}\right)(I_{\mathbb H_t^{N_1}}-yA_2)^{-1}B_2=0,\quad y\in(-\e,\e).
   \]
   Since the pair $(A_2,B_2)$ is controllable we get
   \[
     \sum_{e=1}^Ea_e^\circledast C_2(I_{\mathbb H_t^{N_2}}-x_eA_2)^{-1}
   \]
   and so $0_{\mathbb H_t^{N_1}}$ is sent to $0_{\mathbb H_t^{N_2}}$, and $V$ is well defined.\\
 
 STEP 2: {\sl  The adjoints of the maps $U$ and $V$ with respect to $[\cdot,\cdot]_{\circledast}$ are given by}

\begin{eqnarray}
          U^{\circledast}\left(\sum_{e=1}^E(I_{\mathbb H_t^{N_2}}-y_eA_2^{\circledast})^{-1}C_2^{\circledast}a_e\right)
          &=&\sum_{e=1}^E(I_{\mathbb H_t^{N_1}}-y_eA_1^{\circledast})^{-1}C_1^{\circledast}a_e\\
          V^{\circledast}\left(\sum_{e=1}^E(I_{\mathbb H_t^{N_2}}-y_eA_2)^{-1}B_2a_e\right)&=
            &\sum_{e=1}^E(I_{\mathbb H_t^{N_1}}-y_eA_1)^{-1}B_1a_e
              \label{vcircledast}
        \end{eqnarray}

        These formulas are direct consequences of Proposition \ref{adj-a}.\smallskip

        STEP 3: {\sl We have $UV^{\circledast}=I_{\mathbb H_t^{N_1}}$ and $V^{\circledast}U=I_{\mathbb H_t^N}$ and $N_1=N_2$.}\\

        It follows from \eqref{U-formula} and \eqref{vcircledast} that $V^{\circledast}U=I_{\mathbb H_t^{N_2}}$ and
        $UV^{\circledast}=I_{\mathbb H_t^{N_1}}$. By Corollary \ref{prod-ab-ba}, $N_1=N_2$. \\

        STEP 4:{ \sl Two minimal realizations are similar.}\smallskip

        From \eqref{U-formula} with  $T=1$ and $y_1=0$ we get $UB_1=B_2$ and from \eqref{v-formula} we get $VC_1^{\circledast}=C_2^\circledast$ and so
        $C_1U^{-1}=C_2$ since from Step 3 we have $UV^{\circledast}=I_{\mathbb H_t^{N}}$. Finally \eqref{cab-1}-\eqref{cab-2} imply that
        for every $j,\ell\in\mathbb N_0$
        \[
C_1A_1^jA_1A_1^\ell B_1=C_2A_2^jA_2A_2^\ell B_2.
\]
Using that $UB_1=B_2$ and $C_1=C_2U$ we get
\[
  C_1A_1^j U^{-1}A_2UA_1^kB_1=  C_1A_1^j A_1A_1^kB_1.
\]
Using the observability of $(C_1,A_1)$ and the controllability of $(A_1,B_1)$ we get $U^{-1}A_2U=A_1$.
    \end{proof}

    
    \begin{proposition}
      The realization
      \begin{equation}
        R(x)=D+xC(I_{\mathbb H_t^N} -xA)^{-1}B
      \end{equation}
      is minimal if and only if the pair $(C,A)$ is observable and the pair $(A,B)$ is controllable.
      \end{proposition}

      \begin{proof}
This follows from the fact that two minimal realizations have their main operators of the same size.
      \end{proof}

      We extend the definition of $I$ to matrix-valued functions by
      \begin{equation}
        (I(R))(x)=I(R(x)).
      \end{equation}

      It is readily seen that all properties of $I$ extend. For instance, for matrix-valued rational functions of compatible sizes,
      \begin{equation}
        I(R_1R_2)=I(R_1)I(R_2).
        \label{IR1R2}
        \end{equation}
      \begin{theorem}
        The realization $R(x)=D+xC(I_{\mathbb H_t^N}-xA)^{-1}B$ is minimal if and only if the realization
          \[
(I(R))(x)=I(D)+xI(C)\left(I_{2N}-x I(A)\right)^{-1}I(B)
\]
is minimal.
\label{thm-I-minimal}
\end{theorem}
      \begin{proof}
This follows from the previous proposition, together with Propositions \ref{C-A-I} and \ref{B-A-I}. 
\end{proof}

The proofs of the following two propositions follow the classical complex case, and will be omitted.

      \begin{proposition}
        Assume
        \[
          R(x)=D+xC(I_{\mathbb H_t^N} -xA)^{-1}B
        \]
        be a minimal realization of the $\mathbb H_t^{n\times n}$-valued rational function $R$, and assume $D$ invertible in $\mathbb H_t^{n\times n}$.
        Then, 
        \begin{equation}
          \label{inverse1}
(R(x))^\circledast=D^{\ct}+B^\circledast (I_{\mathbb H_t^N} -xA^\circledast)^{-1}{C^\circledast}
          \end{equation}
          is a minimal realization of $(R(x))^{\circledast}$.
        \end{proposition}

      \begin{proposition}
        Assume
        \begin{equation}
          R(x)=D+xC(I_{\mathbb H_t^N} -xA)^{-1}B
\label{mini-r}
        \end{equation}
        be a minimal realization of the $\mathbb H_t^{n\times n}$-valued rational function $R$, and assume $D$ invertible in $\mathbb H_t^{n\times n}$.
        Then, with $A^\times$ given by \eqref{a-times}, the realization \eqref{inverse11} of $R^{-1}$, i.e.,
\[
R^{-1}(x)=D^{-1}-D^{-1}C(I_{\mathbb H_t^N} -xA^\times)^{-1}BD^{-1}
\]
is a minimal realization of $R^{-1}$.
        \end{proposition}

      \subsection{Minimal factorizations}
      \label{mini-facto-1}
      Let $R$ be a $\in\mathbb H_t^{n\times n}$-valued rational function invertible at the origin, and with minimal realization \eqref{mini-r}, that is of degree $N$.

\begin{definition}
      The factorization $R=R_1R_2$  into a product of two
      $\in\mathbb H_t^{n\times n}$-valued rational functions of respective degrees $N_1$ and $N_2$ is called minimal if $N=N_1+N_2$.
\end{definition}
Note that, already in the complex setting, a rational function of degree $N>1$ may lack non trivial factorizations. (i.e. with both $N_1>0$ and $N_2>0$. In the paper
\cite{MR82e:47024} and the book \cite{bgk1} was given a description of all minimal realizations of a rational function  analytic at $\infty$ (the proof remains the same
if the function is analytic in a neighborhood of the origin).

\begin{remark}
  In view of Theorem \ref{thm-I-minimal} one could study the minimal factorizations of $R$ by first considering the minimal factorizations of $I(R)$ and then restricting to those
  factorizations which have block entries in $\mathcal H_t$. Such an approach will not help to study $\ct$-unitary factorizations (see
  Section \ref{sec-uni} for the latter) since there is no
  simple connection between $I(A^\ct)$ and $I(A)$ for $A\in\mathbb H_t^{n\times n}$ for $t\not=\pm 1$. See Remarks \ref{new-remark} and \ref{remstar}.
    \end{remark}

Before stating and proving the counterpart of this description in the $\mathbb H_t$ setting, namely Theorem \ref{mini-facto}, we need three important definitions.
     
      \begin{definition}
        A matrix $P\in\mathbb H_t^{N\times N}$ is called a projection if $P^2=P$.
      \end{definition}

      \begin{remark}
        We note that, as in the complex case, it holds that $\ker P\cap{\rm ran} P=\left\{0\right\}$
        and
        \begin{equation}
          \label{decomp}
          \mathbb H_t^N=      {\rm ran}\, P \stackrel{\cdot}{[+]} \ker P,
          \end{equation}
where the sum is direct. Furthermore,  if $P$ is a projection so is $SPS^{-1}$ where $P\in\mathbb H_t^{N\times N}$ is invertible in $\mathbb H_t^{N\times N}$.
\end{remark}
      \begin{definition} (see \cite[(1.7) p. 8]{bgk1} for the case of complex numbers)
        Let $R$ be an $\mathbb H_t^{n\times n}$-valued rational function invertible at the origin, and with minimal realization \eqref{mini-r}, and
        let $A^\times$ as in \eqref{a-times}, that is $A^\times=A-BD^{-1}C$. The projection $\pi\in\mathbb H_t^{n\times n}$ is called a supporting
        projection if
        \begin{equation}
          \label{supporting}
          A(\ker \pi)\subset \ker\pi\quad and \quad A^\times ({\rm ran}\,
          \pi )\subset{\rm ran}\,\pi.
          \end{equation}
        \end{definition}

        \begin{lemma}
          Let $S\in\mathbb H_t^{N\times N}$ be a similarity matrix, and assume that $\pi$ is a supporting projection for the realization
          \[
\begin{pmatrix}A&B\\C&D\end{pmatrix}
\]
Then $S\pi S^{-1}$ is a supporting projection for the realization
\[
  \begin{pmatrix}S&0\\ 0&I_{\mathbb H_t^n}\end{pmatrix} \begin{pmatrix}A&B\\C&D\end{pmatrix}\begin{pmatrix}S^{-1}&0\\ 0&I_{\mathbb H_t^n}\end{pmatrix}=\begin{pmatrix}SAS^{-1}&SB\\ CS^{-1}&D\end{pmatrix}
  \]
        \end{lemma}

        \begin{proof}
          The Schur complement of $SAS^{-1}$ in the new realization is
\[
  A^\times_S=SAS^{-1}-SBD^{-1}CS^{-1}=SA^\times S^{-1}, 
\]
and we wish to verify that
\[
  SAS^{-1}(\ker S\pi S^{-1})\subset \ker S\pi S^{-1}\quad{\rm and}\quad
  \underbrace{SA^\times S^{-1}}_{A^\times_S}({\rm ran} S\pi S^{-1})\subset {\rm ran} S\pi S^{-1}.
\]
Since we have
          \[
{\rm ran}\,S\pi S^{-1}=S{\rm  ran}\, \pi\quad {\rm and}\quad \ker S\pi S^{-1}=S\ker\pi,
\]
this follows directly from \eqref{supporting}.
          \end{proof}
        
        For a node $\theta=\begin{pmatrix}A&B\\C&D\end{pmatrix}$ and a supporting projection $\pi$ we define
        \begin{equation}
          {\rm Pr}_\pi (\theta)=\begin{pmatrix}A|_{\ker(\pi)}&(I_{\mathbb H_t^N}-\pi)B\\
              C|_{\ker (\pi)}&D\end{pmatrix}.
          \end{equation}
        
      \begin{theorem}
        \label{mini-facto}
        (see \cite[Theorem 4.8 p. 90]{bgk1} for the complex setting case). Assuming $D=D_1=I_{\mathbb H_t^n}$, there is a one-to-one correspondence
        between minimal factorizations and supporting projections
        \begin{equation}
          \label{div11}
          \theta={\rm Pr_\pi( \theta)}{\rm Pr_{{I_{\mathbb H_t^N}}-\pi}(\theta)}.
        \end{equation}
        More precisely, given a supporting projection $\pi$ the factorization $R=R_1R_2$ with
        \begin{eqnarray}
          R_1(x)&=&I+xC(I_{\mathbb H_t^N}-xA)^{-1}(I-\pi)B,\\
          R_2(x)&=&I+xC\pi(I_{\mathbb H_t^N}-xA)^{-1}B
        \end{eqnarray}
        is minimal and every minimal factorization with factors taking valued $I$ at the origin is of this form.
        \end{theorem}

        \begin{proof}
          We divide the proof in a number of steps, and first assume that we have a minimal factorization.\\

          STEP 1: {\sl Assume $R$ of degree $N$ and $R=R_1R_2$ is a minimal factorization of $R$, where $R_1$ and $R_2$ have minimal realizations as in \eqref{prod_1}, with $D_1=D_2=I_{\mathbb H_t^n}$:
            \[
            R(x)=   \underbrace{ (I_{\mathbb H_t^n}+xC_1(I_{\mathbb H_t^{N_1}}-xA_1)^{-1}B_1)}_{R_1(x)}\underbrace{(I_{\mathbb H_t^n}+
              xC_2(I_{\mathbb H_t^{N_2}}-xA_2)^{-1}B_2)}_{R_2(x)}.
              \]
              Then, \eqref{prod_2}-\eqref{prod_6}
              \begin{eqnarray*}
                  A&=&\begin{pmatrix}A_1&B_1C_2\\
    0_{\mathbb H_t^{N_1}}&A_2\end{pmatrix},\\
  B&=&\begin{pmatrix}B_1\\ B_2\end{pmatrix},\\
  C&=&\begin{pmatrix}C_1& C_2\end{pmatrix},\\
  D&=&I_n,
\end{eqnarray*}
              is a minimal realization of $R$.}\\

Indeed, the above gives a realization of $R=R_1R_2$ by direct computation, as recalled in Section \ref{sec:real}, and the realization is minimal since the factorization is assumed minimal and so $N=N_1+N_2$ is the degree of $R$.\\
            
          STEP 2: {\sl $\pi=\begin{pmatrix}0&0\\0&I_{\mathbb H_t^{N_2}}\end{pmatrix}$ is a supporting projection corresponding to the minimal realization \eqref{prod_2}-\eqref{prod_6} of $R$.}\\

          Indeed,
          \[
            \ker \pi={\rm ran}\,\begin{pmatrix}I_{\mathbb H_t^{N_1}}\\0\end{pmatrix}\quad{\rm and}\quad            {\rm ran}\, \pi={\rm ran}\,
            \begin{pmatrix}0 \\ I_{\mathbb H_t^{N_2}}\end{pmatrix},
          \]
          and \eqref{supporting} are met with $A$ and $A^\times$ as in \eqref{prod_3} and \eqref{atimes}, namely (recall that $D_1=D_2=
          I_{\mathbb H_t^N}$)
          \[
             A=\begin{pmatrix}A_1&B_1C_2\\
             0_{\mathbb H_t^{N_2\times N_1}}&A_2\end{pmatrix}\quad {\rm and}\quad A^\times=\begin{pmatrix}A_1-B_1C_1
             &0_{\mathbb H_t^{N_1\times N_2}}\\    -B_2C_1&A_2-B_2C_2\end{pmatrix}
\]
where the zero matrix in the operator $A$ is of dimensions $N_2\times N_1$ and in the operator $A^\times$ of dimensions $N_1\times N_2$
respectively.\\

We now suppose that we have a supporting projection, and show that it leads to a minimal factorization.\\

STEP 3: {\sl \eqref{supporting} holds.}\\

Consider the block matrix decompositions of $A$ and $A^\times$ along \eqref{decomp}. Condition \eqref{supporting} implies that
\[
A=\begin{pmatrix}A_{11}&A_{12}\\ 0&A_{22}\end{pmatrix}\quad{\rm and}\quad A^\times=\begin{pmatrix}E_{11}&0\\ E_{21}&E_{22}\end{pmatrix}
\]
for some matrices $A_{11},\ldots, E_{22}$.
Set
\[
  B=\begin{pmatrix}B_1\\ B_2\end{pmatrix}\quad{\rm and}\quad C=\begin{pmatrix}C_1&C_2\end{pmatrix}
\]
the condition $A^\times=A-BC$, i.e.
\[
  \begin{pmatrix}E_{11}&0\\ E_{21}&E_{22}\end{pmatrix}=\begin{pmatrix}A_{11}&A_{12}\\ 0&A_{22}\end{pmatrix}-\begin{pmatrix}B_1\\ B_2\end{pmatrix}
  \begin{pmatrix}C_1&C_2\end{pmatrix}
  \]
gives
\[
A_{12}=B_1C_2,
  \]
  so that
  \[
  A=\begin{pmatrix}A_{11}&B_1C_2\\ 0&A_{22}\end{pmatrix}\quad{\rm and}
  \quad A^\times=\begin{pmatrix}A_{11}-B_1C_1&0\\ -B_2C_1&A_{22}-B_2C_2
    \end{pmatrix}.
    \]
    and formulas \eqref{prod_1}-\eqref{prod_0} give
    \[
      R(x)=(I_n+xC_1(I_{N_1}-xA_{11})^{-1}B_1)(I_n+xC_2(I_{N_2}-xA_{22})^{-1}B_2.
      \]
    
      STEP 4: {\sl We have indeed a description of all minimal factorizations of $R$ normalized to be $I$ at the origin.}\smallskip

      It suffices to read the previous arguments backwards.
          \end{proof}
      
          \section{Matrix-rational functions $\ct$-unitary on the imaginary line}
          \setcounter{equation}{0}
\label{sec-uni}
          In the present section and the following ones we consider $\mathbb H_t^{n\times n}$-valued function rational with metric constraints. We first
          consider the case of functions which we call {\it $\ct$-unitary on the imaginary line}. In fact, and following \cite{zbMATH06658818} we replace
          the conditions \eqref{J-unit} (for complex numbers) by the constraints
          \begin{equation}
          R(x)JR(-x)^\ct=J,\quad x\in\mathbb R\cap{\rm Hol}(R),
\label{unit-R}
\end{equation}
where $J\in\mathbb H_t^{n\times n}$ is a signature matrix, that is $J=J^\ct=J^{-1}$.

\subsection{Realization theorem}
Our first result is on the minimal realizations of $\mathbb H_t^{n\times n}$-valued function rational regular at the origin and which satisfy \eqref{unit-R}.

\begin{theorem}
  Let $R$ be a $\mathbb H_t^{n\times n}$-valued rational function regular at the origin, and let $R(x)=D+xC(I_{\mathbb H_t^N}-xA)^{-1}B$ be a minimal realization of $R$. Then $R$ satisfies \eqref{unit-R} if and only if the following conditions hold:\\
$(1)$  $D$ is $\ct$-$J$-unitary, meaning that
\begin{equation}
  \label{lyap-ct-3}
  DJD^\ct= J.
\end{equation}
$(2)$   There exists a $\ct$-Hermitian matrix $H\in\mathbb H_t^{N\times N}$
  such that
  \begin{eqnarray}
    \label{lyap-ct}
    A^\ct H+HA&=&-C^\ct JC\\
    \label{lypa-ct-2}
    B&=&-H^{-1}C^\ct JD.
    \end{eqnarray}
Furthermore, $H$ is uniquely determined by the given minimal realization, and it holds that
\begin{equation}
  \label{JRJR}
      \frac{J-R(x)JR(y)^{\ct}}{x+y}=C(I_{\mathbb H_t^N}-xA)^{-1}H^{-1}(I_{\mathbb H_t^N}-yA^{\ct})^{-1}C^\ct.
    \end{equation}
\label{88-24}
\end{theorem}

  \begin{proof} Rewriting \eqref{unit-R} as
\[
R(x)JR(-x)^\ct J=I_{\mathbb H_t^n},\quad x\in\mathbb R\cap{\rm Hol}(R)
  \]
    we first note that $R(x)$ is invertible for $x\in\mathbb R\cap{\rm Hol}(R)$, with inverse $J(R^{-1}(-x))^\ct J$.
    With $A^\times=A-BD^{-1}C$ and using \eqref{inverse1} for the minimal realization \eqref{inverse11} we have 
    \begin{equation}
      D+xC(I_{\mathbb H_t^N}-xA)^{-1}B=JD^{-\circledast}J+xJD^{-\circledast}B^\circledast(I_{\mathbb H_t^N}+x(A^\times)^\circledast)^{-1}C^\circledast D^{-\circledast}J,
    \end{equation}
    valid for real $x$ where the inverses make sense, where $D^{-\circledast}=(D^{\circledast})^{-1}$ is as in~\eqref{A-ct}. The above is an equality between two minimal realizations of a given matrix-valued rational function and thus there exists a uniquely defined similarity matrix, say $S$, connecting them, that is satisfying the equations
      \begin{eqnarray}
        D&=&JD^{-\ct}J\, ,\\
        A&=&S^{-1}(-A^\times)^\ct S\, ,\\
        \label{eq3!!!}
        C&=&JD^{-\ct}B^\ct S\, ,\\
        SB&=&C^\ct D^{-\ct}J.
      \end{eqnarray}
    Rewriting the second equation as
    \[
      \begin{split}
        SA+A^\ct S&=C^\ct D^{-\ct}B^\ct\\
        &=C^\ct JC \quad ({\rm by}\,\,\, \eqref{eq3!!!})\, ,
      \end{split}
    \]
    we check that these equations can be written in an equivalent way as
    \[
      \begin{split}
      D^\ct JD&=J\, ,\\
      SA+A^\ct S&=C^\ct JC\, ,\\
      SB&=C^\ct D^{-\ct}J\, ,\\
      D^\ct JC&=B^\ct S\, .
\end{split}
      \]
    
      One sees that together with $S$ the matrix $S^\ct$ also satisfies these equations. From the uniqueness of the similarity matrix, we have $S=S^\ct$, which we denote by $-H$. The result follows.\smallskip

      The proof of \eqref{JRJR} is direct, using the Lyapunov equation \eqref{lyap-ct}, and will be omitted.
    \end{proof}

    Following \cite[p. 181]{ag} we define:
    
\begin{definition}
The matrix $H$ in the previous theorem is called the associated $\ct$-symmetric matrix associated to the given minimal realization of $R$.
\end{definition}

As a first example we present the counterpart of a Blaschke factor.

\begin{example} 
\label{bla-bla}
  We take 
\[
A= \alpha\not=-\alpha^{\ct},\quad {\rm and}\quad  J=D=C=1.
\]
Then
\[
  h=\frac{-1}{\alpha+\alpha^\ct}\quad {\rm and}\quad B=\alpha+\alpha^\ct
\]
and the corresponding function $R$ is given by
\begin{equation}
R(x)=(1+x\alpha^\ct)(1-x\alpha)^{-1}.
\end{equation}
When $\alpha$ is not invertible this is a new type of Blaschke factor, different from the classical setting. When $\alpha$ is invertible, setting $\alpha=\beta^{-1}$ we get the counterpart of the more familiar form for the complex setting (see e.g. \cite{Hoffman})
\begin{equation}
b(x)=-\beta \beta^{-\ct}(x+\beta^\ct)(x-\beta)^{-1}.
  \end{equation}
  \end{example}

  \begin{example}
    Let $\alpha,\beta\in\mathbb H_t$ be such that $\alpha+\beta^\ct\not=0$, and set
    \[
H=-J=\begin{pmatrix}0&1\\1&0\end{pmatrix},\quad A=\begin{pmatrix}\alpha&0\\0&\beta\end{pmatrix}\quad{\rm and}\quad  C=\begin{pmatrix}\alpha+\beta^\ct&0\\ 0&1\end{pmatrix}.
\]
Then, \eqref{lyap-ct} is in force,
\[
\begin{split}
  A^\ct H+HA&=\begin{pmatrix}\alpha^\ct&0\\ 0&\beta^\ct\end{pmatrix}\begin{pmatrix}0&1\\1&0\end{pmatrix}+
    \begin{pmatrix}0&1\\1&0\end{pmatrix}\begin{pmatrix}\alpha&0\\ 0&\beta\end{pmatrix}\\
        &=\begin{pmatrix}0&\alpha^\ct+\beta\\
        \alpha+\beta^\ct &0\end{pmatrix}\\
        &=-\begin{pmatrix}\alpha^\ct+\beta&0\\ 0&1\end{pmatrix}\begin{pmatrix} 0&-1\\ -1& 0\end{pmatrix}
          \begin{pmatrix}\alpha+\beta^\ct&0\\ 0&1\end{pmatrix}.
  \end{split}
\]

Taking into account \eqref{lypa-ct-2}, i.e. $B=-H^{-1}C^\ct JD$, we see that the corresponding function $R$ is given by 
\begin{equation}
  \begin{split}
    R(x)&=x\begin{pmatrix}\alpha+\beta^\ct&0\\0&1\end{pmatrix}\begin{pmatrix}(1-\alpha x)^{-1}&0\\0&(1-\beta x)^{-1}\end{pmatrix}
    \begin{pmatrix}0&1\\1&0\end{pmatrix}\begin{pmatrix}\alpha^\ct+\beta&0\\0&1\end{pmatrix}
        \begin{pmatrix}0&1\\1&0\end{pmatrix}\\
          &=\begin{pmatrix}1+x(\alpha+\beta^\ct) (1-\alpha)^{-1}&0\\
          0&1+x(1-\beta x)^{-1}(\alpha^\ct+\beta)\end{pmatrix}\\
                &=\begin{pmatrix}1+(x\alpha-1+1+x\beta^\ct) (1-x\alpha)^{-1}&0\\
          0&1+(1-\beta x)^{-1}(x\alpha^\ct+1-1+x\beta)\end{pmatrix}\\
&=    \begin{pmatrix}(1+x\beta^\ct)(1-x\alpha)^{-1}&0\\  0&(1-x\beta)^{-1}(1+x\alpha^\ct)\end{pmatrix},
  \end{split}
  \label{123rtyui}
\end{equation}
and the $\ct$-$J$-unitarity can also be checked directly:
\[
  \begin{split}
    R(x)JR(-x)^\ct&=\begin{pmatrix}(1+x\beta^\ct)(1-x\alpha)^{-1}&0\\  0&(1-x\beta)^{-1}(1+x\alpha^\ct)\end{pmatrix}\begin{pmatrix}0&1\\1&0\end{pmatrix}\cdot\\
    &\hspace{10mm}\cdot\begin{pmatrix}(1-x\beta^\ct)(1+x\alpha)^{-1}&0\\  0&(1+x\beta)^{-1}(1-x\alpha^\ct)\end{pmatrix}^\ct\\
    &=\begin{pmatrix}0&(1+x\beta^\ct)(1-x\alpha)^{-1}\\ (1-x\beta)^{-1}(1+x\alpha^\ct)&0\end{pmatrix}
    \cdot\\
            &\hspace{10mm}\cdot
    \begin{pmatrix}(1+x\alpha^\ct)^{-1}(1-x\beta)&0\\  0&(1-x\alpha)(1+x\beta^\ct)^{-1}\end{pmatrix}\\
    &=\begin{pmatrix}0&1\\1&0\end{pmatrix}.    
  \end{split}
  \]
See \cite[(7.4) and (7.5) p. 61]{ad3} for a related example in the complex setting. A (non-equivalent) variation is given by
\begin{equation}
  \begin{split}
    R(x)=\begin{pmatrix}(x-\beta^\ct)(x-\alpha)^{-1}&0\\  0&(x-\beta)^{-1}(x-\alpha^\ct)\end{pmatrix}
  \end{split}
  \label{123rtyui2}
        \end{equation}
  \end{example}
  which would be obtained from the counterpart of Theorem \ref{88-24} for minimal realizations centered at infinity. More precisely, replacing $x$ by $1/x$ we obtain from Theorem \ref{88-24} the following result:

  \begin{theorem}
  Let $R$ be a $\mathbb H_t^{n\times n}$-valued rational function regular at infinity, and let $R(x)=D+C(xI_{\mathbb H_t^N}-A)^{-1}B$ be a minimal realization of $R$. Then $R$ satisfies \eqref{unit-R} if and only if the following conditions hold:\\
  $(1)$  $D$ is $\ct$-$J$-unitary.\\
  $(2)$   There exists a $\ct$-symmetric matrix $H\in\mathbb H_t^{N\times N}$
  such that \eqref{lyap-ct} and \eqref{lypa-ct-2} hold.\smallskip

Furthermore, $H$ is uniquely determined by the given minimal realization, and it holds that
\begin{equation}
  \label{JRJRJR}
      \frac{J-R(x)JR(y)^{\ct}}{x+y}=C(xI_{\mathbb H_t^N}-xA)^{-1}H^{-1}(yI_{\mathbb H_t^N}-A^{\ct})^{-1}C^\ct.
    \end{equation}
\label{88-2412}
\end{theorem}

  The next example pertains to the notion of Brune section in classical network theory (see \cite[p. 14]{MR39:5243}), or Blaschke-Potapov factors of the third kind in function
  theory; see \cite{pootapov} for the latter.

  \begin{example}
    We take
    \[
      D=I_{\mathbb H_t^2},\quad
J=\begin{pmatrix}1&0\\0&-1\end{pmatrix}, \quad H\in\mathbb R\setminus \left\{0\right\}, \quad C=\begin{pmatrix}\beta\\ \gamma\end{pmatrix}\quad{\rm and}\quad A=\alpha,
\]
where $\alpha,\beta,\gamma\in\mathbb H_t$ such that
\[
  \alpha=-\alpha^\ct\not=0, \quad {\rm and}\quad\beta\beta^\ct=\gamma\gamma^\ct\not=0.
\]
Then, $C^\ct JC=0$ and the Lyapunov equation is satisfied. The corresponding function $R$ is given by
\begin{equation}
  R(x)=\begin{pmatrix}1&0\\0&1\end{pmatrix}-
  \frac{x}{h}\begin{pmatrix}\beta\\ \gamma\end{pmatrix}(1-\alpha x)^{-1}\begin{pmatrix}\beta^\ct &-\gamma^\ct\end{pmatrix}.
  \end{equation}
    \end{example}
\subsection{$\ct$-unitary minimal factorizations}
\label{c45}
We first remark:

\begin{lemma}
  In the above notation, it holds that
  \begin{equation}
    \label{indefinite}
H^{-1}A^{\ct}H=-A^\times,
\end{equation}
and in particular $A^\times = -A^\Delta$, where $A^\Delta$ is the adjoint of $A$ in the form \eqref{[]H}.
  \end{lemma}

  \begin{proof} We write
    
    \[
      \begin{split}
      A-BD^{-1}C&=A+H^{-1}C^\ct JD D^{-1}C,\\
      &=A+H^{-1}C^\ct JC,\\
      &=-H^{-1}A^\ct H.
      \end{split}
    \]
  \end{proof}

  The following result is the counterpart of \cite[Theorem 2.6 p 187]{ag} in the present setting.

  \begin{theorem}
    \label{facto-1234567}
    Let $R(x)=D+xC(I_{\mathbb H_t^N}-xA)^{-1}B$ satisfying \eqref{unit-R}, with associated $\ct$-symmetric matrix $H$. Up to $J$-unitary factorizations of $D$ there is a
      one-to-one correspondence between $J$-unitary minimal factorizations of $R$ and $A$-invariant subspaces which are non-degenerate in the $H$-metric.
    \end{theorem}

    \begin{proof} Without loss of generality we can, and will, assume that $D=I_{\mathbb H_t^n}$ and that all factors have also this normalization. By Theorem \ref{88-24} there
      exists a uniquely defined invertible $\ct$-symmetric matrix $H\in\mathbb H_t^{N\times N}$ such that \eqref{lyap-ct} and \eqref{lypa-ct-2} hold, that is (with $D=I_n$)
      \begin{eqnarray}
        \nonumber
    A^\ct H+HA&=&-C^\ct JC,\\
    \label{lypa-ct-2-2}
    HB&=&C^\ct J.
      \end{eqnarray}
      Assume now $R=R_1R_2$ denote a non-trivial $J$-unitary minimal factorization. By Theorem \ref{mini-facto}, there exists a supporting projection $\pi$ such that
      \eqref{supporting} holds. In the present case we have $A^\times=A^\Delta$. By Theorem \ref{lemma-preumss}, the supporting projection is $H$-self-adjoint, i.e.
      \[
\pi^\ct H=H\pi^\ct
        \]
        (see \cite[(2.36) p. 189]{ag} for the classical counterpart) . But
        \begin{equation}
          \label{pi-1}
          \pi=\begin{pmatrix}0&0\\0&I_t^{N_2}\end{pmatrix}
\end{equation}
          and it follows that $H$ is block diagonal
        \[
          H=\begin{pmatrix} H_{11}&0\\0&H_{22}\end{pmatrix}
        \]
        where $H_{11}\in\mathbb H_t^{N_1\times N_1}$. Multiplying both sides \eqref{lyap-ct} on the left and on the right
        by $\pi$ and multiplying both sides of \eqref{lypa-ct-2-2} on the left by $\pi$ we obtain
        \begin{eqnarray}
          A_{11}^\ct H_{11}+H_{11}A_{11}&=&C_1^\ct J C_1\\
          -H_{11}B_1&=&C_1^\ct J,
                      \end{eqnarray}
                      and it follows that the factor  $R_1$ is $J$-unitary. The factor $R_2$ is then also $J$-unitary.
    \end{proof}

    Not every $J$-unitary function will admit non-trivial $J$-unitary factorizations, although it may admit non-trivial minimal factorizations. An example is provided
    by the function \eqref{123rtyui}.
    \section{$\ct$-rational functions $\ct$-unitary on the unit circle}
    \setcounter{equation}{0}
  As explained in Remark \ref{remark-r-t-1}, in the complex case, condition \eqref{J-unitary-123} can be rewritten as
    \[
      R(z)JR(1/\overline{z})^*=J
    \]
    at those points of the complex plane where the functions are defined. In the setting of $\mathbb H_t$, and following \cite{zbMATH06658818}
    where the quaternionic setting is considered, we define $\ct$-unitarity on the circle by the equation
      \begin{equation}
        \label{J-TTT}
        R(x)JR((1/x))^\ct =J
      \end{equation}
    at those points where the functions are defined.
    \subsection{Realization theorem}
    Already in the complex setting the circle case is a bit more involved and requires  conditions of invertibility. In the complex setting these conditions can be relaxed using the theory of finite dimensional de Branges Rovnyak spaces Pontryagin spaces; see \cite{ad3,ad-laa3}.
    We will consider this aspect in the setting of the rings $\mathbb H_t$ in a separate publication.

\begin{theorem}
      Let $R$ be a $\mathbb H_t^{n\times n}$-valued rational function regular at the origin, and let $R(x)=D+xC(I_{\mathbb H_t^N}-xA)^{-1}B$ be a minimal realization of $R$. Assume moreover that $R$ is analytic and invertible at infinity. Then $R$ satisfies \eqref{J-TTT}
      if and only if there exists a $\ct$-symmetric matrix $H\in\mathbb H_t^{N\times N}$ such that
      \begin{equation}
        \begin{pmatrix}          A&B\\ C&D\end{pmatrix}        \begin{pmatrix}          H^{-1}&0\\ 0&J\end{pmatrix}
        \begin{pmatrix}          A&B\\ C&D\end{pmatrix}^\ct=      
        \begin{pmatrix}
          H^{-1}&0\\ 0&J\end{pmatrix}.
        \label{unit-circle-4}
 \end{equation}
 \label{unit-circle}
\end{theorem}

    \begin{proof} We follow the proof of \cite[Theorem 3.1 p. 197]{ag} and proceed in a number of steps.\\

      STEP 1: {\sl The main operator $A$ in the given minimal realization of $R$ is invertible.}\smallskip

      By Theorem \ref{thm-I-minimal} the realization
      \[
        (I(R))(x)=I(D)+xI(C)\left(I_{2N}-x I(A)\right)^{-1}I(B)
      \]
      is minimal. Using the definition of convergence and the map $I$ ( see Section \ref{conv-123}) we see that $I(R)(\infty)$ exists and is invertible. It follows that $I(R)(1/x)$ is regular and invertible at the origin, and by e.g. \cite[Corollary 2.7 p. 59]{bgk1}, $I(A)$ is invertible. By Corollary \ref{ab-ba}
      it follows that $A$ is invertible.\\

      STEP 2: {\sl 
A minimal realization of $JR(1/x)^\ct J$ is given by   
      \[
        JR(1/x)^\ct J=J(D^\ct-B^\ct (A^\ct)^{-1}C^\ct)J-xB^\ct (A^\ct)^{-1}(I_{\mathbb H_t^N}-x(A^\ct)^{-1})^{-1}(A^{\ct})^{-1}C^\ct J.
      \]
    }
We recall (see Corollary \ref{ab-ba}) that $A$ is invertible if and only if $A^\ct$ is invertible. We can thus write
\[
  \begin{split}
    JR(1/x)^\ct J&=JD^\ct J+B^\ct (xI_{\mathbb H_t^N}-A^\ct)^{-1}C^\ct J\\
    &=JD^\ct J-B^\ct (A^{\ct})^{-1} (I_{\mathbb H_t^N}-x(A^\ct)^{-1})^{-1}C^\ct J\\
    &=JD^\ct J-B^\ct (A^{\ct})^{-1}C^\ct J +JB^\ct (A^{\ct})^{-1}\left\{I_{\mathbb H_t^N}-(I_{\mathbb H_t^N}-x(A^\ct)^{-1})^{-1}\right\}C^\ct J\\
    &=JD^\ct J-B^\ct (A^{\ct})^{-1}C^\ct J -xJB^\ct (A^{\ct})^{-1}(I_{\mathbb H_t^N}-x(A^\ct)^{-1})^{-1}(A^\ct)^{-1}C^\ct J,
  \end{split}
\]
and hence the result.\smallskip

STEP 3: {\sl There exists an invertible matrix $S\in\mathbb H_t^{N\times N}$ such that
  \begin{equation}
\label{S-simi}
    \begin{split}
    \begin{pmatrix}
      A-BD^{-1}C&BD^{-1}\\ -D^{-1}C&D^{-1}\end{pmatrix}\begin{pmatrix}S&0\\ 0&I_{\mathbb H_t^n}\end{pmatrix}=\\
    &\hspace{-3cm}
    =\begin{pmatrix}S&0\\ 0&I_{\mathbb H_t^n}\end{pmatrix}\begin{pmatrix} (A^\ct)^{-1}&C^\ct J\\ -JB^\ct (A^\ct)^{-1}&J(D^\ct-B^\ct (A^\ct)^{-1}C^\ct)J\end{pmatrix}.
    \end{split}
  \end{equation}
}

We rewrite \eqref{J-TTT} as

\begin{equation}
  \label{1987}
  R(x)^{-1}=JR(1/x)^\ct J.
  \end{equation}
  A minimal realization of $R^{-1}$ in terms of the original minimal realization is given in \eqref{inverse1}, and a minimal realization of $JR(1/x)^\ct J$ has been given in the previous step.
  The uniqueness of the minimal realization up to a similarity matrix applied to \eqref{1987} allows to conclude.\smallskip

STEP 4: {\sl $S=S^\ct$ and \eqref{unit-circle-4} is in force with $S=-H^{-1}$.}\smallskip

Indeed, slightly reorganizing \eqref{J-TTT} we see that the latter is equivalent to:
\begin{eqnarray}  
    \label{labek1}
S-ASA^\ct&=&-BD^{-1}CSA^\ct,\\
  \label{labek2}
  CSA^\ct&=&DJB^\ct,\\
    \label{labek3}
  B&=&SC^\ct JD,\\
    \label{labek4}
  J&=&DJD^\ct-DJB^\ct (A^\ct)^{-1}C^\ct,
\end{eqnarray}
or, equivalently,
\begin{eqnarray}  
    \label{labek11}
S-ASA^\ct&=&-BJB^\ct,\\
  \label{labek21}
  CSA^\ct&=&DJB^\ct,\\
    \label{labek31}
  B&=&SC^\ct JD,\\
    \label{labek41}
  J&=&DJD^\ct-CSC^\ct,
\end{eqnarray}
that is,
\begin{equation}
        \begin{pmatrix}          A&B\\ C&D\end{pmatrix}        \begin{pmatrix}          -S&0\\ 0&J\end{pmatrix}
        \begin{pmatrix}          A&B\\ C&D\end{pmatrix}^\ct=      
   \begin{pmatrix}         - S&0\\ 0&J\end{pmatrix}.
 \label{unit-circle-1}
  \end{equation}

  Taking $\ct$ of both sides of this equality we see that $S^\ct$ also solves \eqref{unit-circle-1}. By uniqueness of the similarity matrix we have $S=S^\ct$ and hence, comparing with \eqref{unit-circle-4}, $S=-H^{-1}$.  
    \end{proof}

    We note that \eqref{unit-circle-4} is equivalent to

    \begin{equation}
        \begin{pmatrix}          A&B\\ C&D\end{pmatrix}^{\ct}       \begin{pmatrix}          H&0\\ 0&J\end{pmatrix}
        \begin{pmatrix}          A&B\\ C&D\end{pmatrix}=      
   \begin{pmatrix}          H&0\\ 0&J\end{pmatrix}.
 \label{unit-circle-8}
 \end{equation}

 \begin{corollary}
   \begin{equation}
\frac{J-R(x)JR(y)^\ct}{1-xy}=C(I_{\mathbb H_t^N}-xA)^{-1}H^{-1}(I_{\mathbb H_t^N}-yA^\ct)^{-1}C^\ct.
     \end{equation}
   \end{corollary}

   \begin{proof} Using \eqref{labek11}-\eqref{labek41} to go from the third line to the following lines we write
(we replace $S$ by $-H^{-1}$ at the end of the computations)
     \[
\begin{split}
  R(x)JR(y)^\ct-J&=(D+xC(I_{\mathbb H_t^N}-xA)^{-1}B)J(D^\ct+yB^\ct (I_{\mathbb H_t^N}-yA^\ct)^{-1}C^\ct)\\
  &=DJD^\ct+xC(I_{\mathbb H_t^N}-xA)^{-1}BJD^\ct+yDJB^\ct(I_{\mathbb H_t^N}-yA^\ct)^{-1}C^\ct+\\
  &\hspace{5mm} +xC(I_{\mathbb H_t^N}-xA)^{-1}BJB^\ct(I_{\mathbb H_t^N}-yA^\ct)^{-1}C^\ct-J\\
  &=CSC^*+xC(I_{\mathbb H_t^N}-xA)^{-1}ASC^\ct+yCSA^\ct(I_{\mathbb H_t^N}-yA^\ct)^{-1}C^\ct+\\
  &\hspace{5mm} +xC(I_{\mathbb H_t^N}-xA)^{-1}(ASA^\ct-S)(I_{\mathbb H_t^N}-yA^\ct)^{-1}C^\ct\\
  &=C(I_{\mathbb H_t^N}-xA)^{-1}\left\{\diamondsuit\right\}(I_{\mathbb H_t^N}-yA^\ct)^{-1}C^\ct,
\end{split}
\]
with
\[
  \begin{split}
    \diamondsuit&=(I_{\mathbb H_t^N}-xA)S(I_{\mathbb H_t^N}-yA^\ct)+xAS(I_{\mathbb H_t^N}-yA^\ct)+\\
    &\hspace{5mm}+y(I_{\mathbb H_t^N}-xA)SA^\ct+xy(ASA^\ct-S)\\
    &=(1-xy)S.
  \end{split}
  \]
\end{proof}

As a first example we give the Blaschke factor, defined in \cite[Definition 4.7]{acv1} and earlier in \cite{MR2872478,MR3192300,zbMATH06658818} in the setting of quaternions. For
$\alpha\in\mathbb H_t$ such that $\alpha\alpha^\ct\not=1$ we define
\begin{equation}
  b_\alpha(x)=(x-\alpha)(1-x\alpha^\ct)^{-1}.
\end{equation}
Rewriting
\[
  \begin{split}
b_\alpha(x)&=-\alpha+x(1-\alpha\alpha^\ct)(1-\alpha^\ct)^{-1}
    \end{split}
  \]
  we get, with
  \begin{equation}
    \begin{split}
      A&=\alpha^\ct, \\
      B&=1,\\
      C&=1-\alpha\alpha^\ct,\\
      D&=-\alpha,
      \end{split}
    \end{equation}
    and $h=1-\alpha\alpha^\ct$ (recall that $\alpha\alpha^\ct$ and hence $h$, is real; see Lemma \ref{unit-p-0}), 
    \[
      \begin{pmatrix}A&B\\ C&D\end{pmatrix}^\ct\begin{pmatrix}h&0\\0&1\end{pmatrix}
      \begin{pmatrix}A&B\\ C&D\end{pmatrix}^\ct=\begin{pmatrix}h&0\\0&1\end{pmatrix}.
    \]
When $h>0$ we have the counterpart of a classical Blaschke factor while $h<0$ correspond to the inverse of a classical Blaschke factor.\\

A finite product (with respect to the real variable $x$ or $\star$-product for the variable $q$) will give also a unitary function. As in the quaternionic case, and as shown in \cite{acv1}, the case $b_\alpha\star b_{\alpha^\ct}$ is of special interest. One has
(see \cite[Lemma 4.9]{acv1})
\begin{equation}
(b_\alpha\star b_{\alpha^\ct})(q)=\frac{q^2-2q{\rm Re}\alpha+\det \alpha}{q^2\det\alpha-2q{\rm Re}\alpha+1},\quad q\in\mathbb H_t,
  \end{equation}
  where the order of the factors is not important since the coefficients are real.\\
  
As for the line case, another (non-equivalent) family of Blaschke factors is given by $(1-x\beta)(x-\beta^\ct)^{-1}$. In the present case, this is just
$b_{\beta^\ct}^{-1}$.\\

The following example is also taken from our previous paper \cite{acv1}.

\begin{example}
  Let
                           \begin{eqnarray}
\label{theta-!!!}
                             \Theta(q)&=&D+qC\star(I_{\mathbb H_1^N}-qA)^{-\star}B
                             \end{eqnarray}
                           where 
         
         \begin{eqnarray}
           A&=&{\rm diag}(\alpha_1^{\ct},\ldots, \alpha_N^{\ct})
           \label{t-1}\\
                B&=&(I_{\mathbb H_t^N}-A){\mathbf G_{\ct}}^{-1}(I_{\mathbb H_t^N}-A^{\ct})^{-1}C^{\ct}\\
                     C&=&\begin{pmatrix}1&1&\cdots &1\end{pmatrix}\\
           D&=&1-C{\mathbf G_{\ct}}^{-1}(I_{\mathbb H_t^N}-A^{\ct})^{-1}C^{\ct}.
                \label{t-4}
         \end{eqnarray}

         The matrix  ${\mathbf G_{\ct}}$ is the unique solution of the Stein equation
  \begin{equation}
    {\mathbf G_{\ct}}-A^{\ct}{\mathbf G_{\ct}} A=C^{\ct}C,
    \label{stein-eq}
  \end{equation}
  and is given by the formula
  \begin{equation}
    \label{formula-G}
    {\mathbf G_{\ct}}=\sum_{n=0}^\infty A^{\ct n}C^\ct CA^n.
  \end{equation}
  Furthermore,
  
  \begin{eqnarray}
  A^\ct{\mathbf G_{\ct}} A+C^\ct C&=&{\mathbf G_{\ct}},
                              \label{t-9}\\
  \label{t-890}
  B^{\ct}{\mathbf G_{\ct}} A+ D^\ct C&=&0,\\
  B^{\ct}{\mathbf G_{\ct}} B +D^\ct D&=&I,
                                 \label{t-900}
    \end{eqnarray}

    i.e. \eqref{unit-circle-8} is in force with $H=\mathbf G_\ct$ and $J=1$.
  \end{example}

  We note that:

  \begin{proposition}
    The function can be rewritten as
    \begin{equation}
      \label{bB-1234}
\Theta(x)=1-(1-x)C(I_{\mathbb H_t^N}-xA)^{-1}G_\ct^{-1}(I_{\mathbb H_t^N}-A^\ct)^{-1}C^\ct.
  \end{equation}
  Furthermore, when $N=1$, we have (and we write $b$ rather than $\Theta$ and $\alpha$ instead of $A$)
  \begin{equation}
    \label{b123}
    b(x)=(1-x\alpha)^{-1}(x-\alpha^\ct)(1-\alpha)(1-\alpha^\ct)^{-1}.
    \end{equation}
    \end{proposition}

    \begin{proof}
      To prove \eqref{bB-1234} we write
      \[
        \begin{split}
          \Theta(x)&=1-CG_\ct^{-1}(I_{\mathbb H_t^N}-A^{\ct})^{-1}C+xC(I_{\mathbb H_t^N}-xA)^{-1}G_\ct^{-1}(I_{\mathbb H_t^N}-A^\ct)^{-1}C^\ct\\
          &=1+C(I_{\mathbb H_t^N}-xA)^{-1}\{\underbrace{x(I_{\mathbb H_t^N}-A)-I_{\mathbb H_t^N}+x}_{=(x-1)I_{\mathbb H_t^N}}\}
          G_\ct^{-1}(I_{\mathbb H_t^N}-A^\ct)^{-1}C^\ct.
        \end{split}
        \]
      
      We now prove \eqref{b123} and write, using \eqref{bB-1234},
      \[
        \begin{split}
          b(x)&=1+(1-x\alpha)^{-1}(x-1)G_\ct^{-1}(1-\alpha^\ct)^{-1}\\
          &=(1-x\alpha)^{-1}\left((1-x\alpha)G_\ct(1-\alpha^\ct)+x-1\right)G_\ct^{-1}(1-\alpha^\ct)^{-1}\\
          &=(1-x\alpha)^{-1}\left((1-x\alpha)(1-\alpha^\ct)+G_\ct^{-1}( x-1)\right)(1-\alpha^\ct)^{-1}\\
                    &=(1-x\alpha)^{-1}\left(\underbrace{(1-x\alpha)(1-\alpha^\ct)+(1-\alpha^\ct\alpha)( x-1)}_{=(x-\alpha^\ct)(1-\alpha)}\right)(1-\alpha^\ct)^{-1},
          \end{split}
        \]
        where we have used that, for $N=1$ and $\alpha_1=\alpha$, equation \eqref{t-9} gives
        \[
          G_\ct=(1-\alpha^\ct\alpha)^{-1}.
          \]
      \end{proof}
\subsection{$\ct$-minimal factorizations}
\label{c57}

Equation \eqref{labek1} (or the upper left block equality in \eqref{S-simi} gives
\begin{equation}
  A^\times A^\Delta=I_{\mathbb H_t^N},
\end{equation}
where we recall that $A^\times=A-BD^{-1}C$ and that $A^\Delta=H^{-1}A^\ct H$ denotes the adjoint of $A$ in the
$[\cdot,\cdot]_H$ symmetric form; see \eqref{Adelta}.\smallskip

The counterpart of Theorem \ref{facto-1234567} is now (note that in the statement, the claim {\sl Up to $J$-unitary factorizations of $D$} is replaced by
{\sl Up to a $J$-unitary factor}). Note that we need the hypothesis that $(I_{\mathbb H_t^N}-A)$ (or $(I_{\mathbb H_t^N}+A)$) is invertible. We use some formulas from  \cite{ag}, with the adjoint replaced by $\ct$.

\begin{theorem}
  Let $R(x)=D+xC(I_{\mathbb H_t^N}-xA)^{-1}B$ satisfying \eqref{unit-R}, with associated $\ct$-symmetric matrix $H$. Assume both $A$ and $I_{\mathbb H_t^N}-A$ invertible.
  Up to a $J$-unitary factor there is a     one-to-one correspondence between $J$-unitary minimal factorizations of $R$ and $A$-invariant subspaces which are non-degenerate
  in the $H$-metric.
  For a given orthogonal supporting projection corresponding to the decomposition \eqref{prod_3}-\eqref{prod_6}, the left $J$-unitary factor is given by
\begin{equation}
  \label{R1-ag}
R_1(x)=(I_n+C_1H_{11}^{-1}(I_{\mathbb H_t^{N_1}}-A_1)^{-1})H_{11}^{-1}(A_{11}^{-1})^\ct C_1^\ct J)(I_n-CH^{-1}(I_{\mathbb H_t^{N_1}}-A_1^\ct)^{-1}C_1^\ct J).
    \end{equation}
\end{theorem}
    \begin{proof}
      The proof goes along the lines of the proof of Theorem \ref{facto-1234567}, using formulas from \cite[(3.20)-(3-21)]{ag}.
Consider \eqref{unit-circle-4}
with $A,B,C$ and $D$ given by \eqref{prod_3}-\eqref{prod_6}. Then $\pi$ is still given by \eqref{pi-1}.  Since it is $H$ $\ct$-symmetric it follows that $H$ is block diagonal (as in the proof of Theorem \ref{facto-1234567}).
Thus it holds that
\begin{equation}
  \label{h112233}
H_{11}-A_1^\ct H_{11}A_1=C_1^\ct JC_1.
  \end{equation}

  Rewrite $R_1$ as
  \[
R_1(x)=\mathscr D+x\mathscr C(I_{\mathbb H_t^{N_1}}-x\mathscr A)^{-1}\mathscr B
    \]
    with
\begin{equation}
  \mathscr A=A_{11},
  \end{equation}
and
\begin{eqnarray}
  \label{B-D}
  \mathscr B&=&-H_{11}^{-1}(A_{11}^{-1})^\ct C_1^\ct J(I_{\mathbb H_t^{n}}-CH^{-1}(\overbrace{I_{\mathbb H_t^{N_1}}-A_1^\ct)^{-1}C_1^\ct J}^{\mathscr D}),\\
  \mathscr C&=&C_1H_{11},\\
  \mathscr D&=&I_{\mathbb H_t^n}-C_1H_{11}^{-1}(I_{\mathbb H_t^{N_1}}-A_1^\ct)^{-1}C_1^\ct J.
\end{eqnarray}

We want to prove that
\begin{equation}
  \begin{pmatrix}\mathscr A&\mathscr B\\ \mathscr C&\mathscr D\end{pmatrix}^\ct\begin{pmatrix}H_{11}&0\\0&J\end{pmatrix}
  \begin{pmatrix}\mathscr A&\mathscr B\\ \mathscr C&\mathscr D\end{pmatrix}=
  \begin{pmatrix}H_{11}&0\\0&J\end{pmatrix}.
\end{equation}
The equality in the upper left block is just \eqref{h112233}, while the equality in the upper right block, which reads
\[
{\mathscr A} H_{11}\mathscr D+\mathscr C J\mathscr D=0
  \]
is just \eqref{B-D}.\smallskip

The formula \eqref{R1-ag} follows then from \cite[(3.20)-(3-21)]{ag}. It remains thus to check the equality for the lower right block, that is
\[
  \mathscr B^\ct H_{11}\mathscr B+\mathscr D^\ct J\mathscr D=J,
\]
which, provided we prove that  that $\mathscr D$ is invertible, we will rewrite as 
\begin{equation}
  \label{new-eq}
\mathscr D^{-\ct}\mathscr B H_{11}\mathscr B\mathscr D^{-1}+J=\mathscr D^{-\ct}\mathscr D^{-1}.
  \end{equation}

  STEP : {\sl $\mathscr D$ is invertible.}\smallskip

  Using the classical formula $(1-ab)^{-1}=1+a(1-ba)^{-1}b$, where $a,b$ are matrices of compatible entries in some ring and the $1$ denote appropriate identity matrices we have:

  \[
    \begin{split}
      \mathscr D^{-1}&=(I_{\mathbb H_t^n}-CH^{-1}(I_{\mathbb H_t^N}-A^\ct)^{-1}C^\ct J)^{-1}\\
      &=I_{\mathbb H_t^n}+CH^{-1}(I_{\mathbb H_t^N}-(I_{\mathbb H_t^N}-A^\ct)^{-1}C^\ct JC H^{-1})^{-1}(I_{\mathbb H_t^N}-A^\ct)^{-1}C^\ct J\\
      &=I_{\mathbb H_t^n}+C(H-A^\ct H-C^\ct JC)^{-1}C^\ct J\\
      &=I_{\mathbb H_t^n}+C(A^\ct HA-A^\ct H)^{-1}C^\ct J\\
      &=I_{\mathbb H_t^n}+C(A-I_{\mathbb H_t^N})^{-1}H^{-1}A^{-\ct}C^\ct J.
%
%
%
%
  %
    \end{split}
  \]

  STEP : {\sl We compute $\mathscr D^{-\ct}J\mathscr D^{-1}-J$.}\smallskip

  We have
  
  \[
    \begin{split}
      \mathscr D^{-\ct}J\mathscr D^{-1}-J
      &=JCA^{-1}H^{-1}(A^\ct-I_{\mathbb H_t^N})^{-1}\left\{\diamondsuit\right\}(A-I_{\mathbb H_t^N})^{-1}H^{-1}A^{-\ct}C^\ct J
\end{split}
\]
with
\[
  \begin{split}
    \diamondsuit&=A^\ct H(A-I_{\mathbb H_t^N})+(A^\ct -I_{\mathbb H_t^N})HA+C^\ct JC\\
    &=A^\ct H A+H-A^\ct H-HA\\
    &=(A^\ct-I_{\mathbb H_t^N})H(A-I_{\mathbb H_t^N})
  \end{split}
\]
and so
\begin{equation}
  \label{eq-67}
  \mathscr D^{-\ct}J \mathscr D^{-1}-J=JCA^{-1}H^{-1}A^{-\ct}C^\ct J.
\end{equation}

  STEP : {\sl \eqref{new-eq} holds.}\smallskip

 Since $\mathscr D^{-1}\mathscr B=H_{11}^{-1}A_1^{-\ct}C_1^\ct J$, we can rewrite \eqref{eq-67} as
\[
  \mathscr D^{-\ct}J \mathscr D^{-1}-J\mathscr D^{-\ct}\mathscr B^\ct H\mathscr B\mathscr D^{-1}
\]
and so

\begin{equation}
    J=\mathscr D^\ct J \mathscr D+\mathscr B^\ct H \mathscr B.
    \label{B-D}
    \end{equation}
      \end{proof}
      \section{$\ct$-rational function $\ct$-anti-symmetric on the imaginary line}
\setcounter{equation}{0}
In this section we characterize minimal realizations and minimal additive decompositions of rational  matrix-valued functions satisfying
\begin{equation}
  \label{cara-sym}
  \phi(x)=-(\phi(-x))^\ct
  \end{equation}
at those real points where $\phi$ is defined.
\subsection{$\ct$-minimal realizations}
\label{phi--phi-ct}

  \begin{theorem}
 The $\mathbb H_t^{n\times n}$-valued rational function $\phi$ with minimal realization $\phi(x)=D+xC(I_{\mathbf H_t^N}-xA)^{-1}B$ satisfies \eqref{cara-sym} at the real points where it is defined if and only if there exists a $\ct$-symmetric matrix $H\in\mathbb H_t^{N\times N}$ such that
 \begin{eqnarray}
   \label{sna1}
      D+D^\ct&=&0,\\
      A^\ct H+HA&=&0,\\
   B&=&H^{-1}C^\ct.
        \label{sna4}
    \end{eqnarray}
    Furthermore, $H$ is uniquely determined by the realization and it holds that
    \begin{equation}
      \phi(x)=D+xC(I_{\mathbb H_t^N}-xA)^{-1}H^{-1}C^\ct
      \label{phi-real}
    \end{equation}
    and
\begin{equation}
  \label{rk-phi}
\frac{\phi(x)+(\phi(y))^\ct}{x+y}=C(I_{\mathbb H_t^N}-xA)^{-1}H^{-1}(I_{\mathbb H_t^N}-yA^{\ct})^{-1}C^\ct.
      \end{equation}
    \end{theorem}
\begin{proof}
  We set, see \eqref{JJJJ},
  \begin{equation}
    \label{J-T}
    J=\begin{pmatrix}0&I_{\mathbb H_t^n}\\ I_{\mathbb H_t^n}&0\end{pmatrix}\quad{\rm and}\quad
    T(x)=\begin{pmatrix}I_{\mathbb H_t^n}&\phi(x)\\0&I_{\mathbb H_t^n}\end{pmatrix}
  \end{equation}
  and divide the proof into a number of steps. The proofs of Steps 1 and 2 are clear, and will be omitted.\\

  STEP 1: {\sl \eqref{cara-sym} holds if and only if $T(x)JT(-x)^{\ct}=J$, with $T(x)$ and $J$ as in \eqref{J-T}.}\\

  STEP 2: {\sl Let $\phi(x)=D+xC(I_{\mathbb H_t^N}-xA)^{-1}B$ be a minimal realization of $\phi$. Then,
    \begin{equation}
      \label{Tx}
T(x)=\mathscr D+x\mathscr C(I_{\mathbb H_t^N}-xA)^{-1}\mathscr B
\end{equation}
is a minimal realization of $T(x)$, where  
\begin{equation}
  \label{Real-T}
  \begin{split}
    \mathscr B&=\begin{pmatrix}0& B\end{pmatrix}\\
    \mathscr C&=\begin{pmatrix}C\\ 0\end{pmatrix}\\
\mathscr D&=\begin{pmatrix}I_{\mathbb H_t^n}& D\\ 0&I_{\mathbb H_t^n} \end{pmatrix}.
    \end{split}
  \end{equation}
}

STEP 3: {\sl We apply Theorem \ref{88-24} to the above realization.}\smallskip

We replace $H$ by $-H$. The state operator $A$ does not change. With $J$ as in \eqref{J-T} we have $\mathscr C^*J\mathscr C=0$ and the Lyapunov equation \eqref{lyap-ct} becomes 
\begin{equation}
A^\ct H+HA=0
\label{ActHHA}
\end{equation}
Furthermore equation \eqref{lyap-ct-3} becomes
\[
  D+D^\ct=0
\]
and \eqref{lypa-ct-2} becomes
\[
\mathscr B=H^{-1}DJ\mathscr C^\ct J,
\]
that is
\begin{equation}
  \label{456}
B=H^{-1}C^\ct.
\end{equation}

STEP 4: {\sl We prove \eqref{phi-real} and  \eqref{rk-phi}.}\smallskip

\eqref{phi-real} follows directly from \eqref{456}. For \eqref{rk-phi} we write:
\[
  \begin{split}
    \phi(x)+\phi(y)^\ct&=D+D^\ct+xC(I_{\mathbb H_t^N}-xA)^{-1}H^{-1}C^\ct+yCH^{-1}(I_{\mathbb H_t^N}-yA)^{-\ct}C^\ct\\
    &\hspace{-2cm}=C(I_{\mathbb H_t^N}-xA)^{-1}H^{-1}\left\{
      x(I_{\mathbb H_t^N}-yA)^{\ct}H+yH(I_{\mathbb H_t^N}-xA)   \right\}H^{-1}(I_{\mathbb H_t^N}-yA)^{-\ct}C^\ct\\
    &\hspace{-2cm}=(x+y)C(I_{\mathbb H_t^N}-yA)^{-1}H^{-1}(I_{\mathbb H_t^N}-yA)^{-\ct}C^\ct,
  \end{split}
  \]
and hence the result.
\end{proof}

\subsection{$\ct$-unitary minimal additive decompositions}
\label{c56}
\begin{theorem}
  \label{thm-line-decomp}
  Let $\phi$ be a $\mathbb H_t^{n\times}$-valued rational function analytic at the origin, and satisfying \eqref{cara-sym},  and let 
  $\phi(x)=D+xC(I_{\mathbb H_t^n}-xA)^{-1}B$ be a minimal realization of $\phi$, with associated $\ct$-symmetric matrix $H$. Then, up to an additive
  $\ct$-anti-symmetric decomposition of $D-D^\ct$ into two $\ct$-anti-symmetric matrices,
\[
  D-D^\ct=M_1+M_2, \quad where\quad M_j+M_j^\ct=0,\,\, j=1,2,
\]
  there is a one-to-one correspondence between $\ct$-unitary minimal additive decompositions of $\phi$ and $H$ non-degenerate $A$-invariant subspaces
    of $\mathbb H_t^N$.
\end{theorem}

\begin{proof}
  Such additive decompositions correspond to minimal factorizations of the function $T$ defined by \eqref{Tx}. The operator $A^\times $ for the minimal
  realization \eqref{Real-T} of $T(x)$ is equal to
  \begin{equation}
    A^\times=A-\mathscr B\mathscr D^{-1}\mathscr C=A,
  \end{equation}
  and equation \eqref{ActHHA}, rewritten as
\[
A=-H^{-1}A^\ct A=-A^\Delta
\]
  expresses that $A$ is $\ct$ anti self-adjoint with respect to the form \eqref{[]H} defined by $H$ (recall that the corresponding adjoint $A^\Delta$ is given by \eqref{Adelta}). To conclude it suffices to apply
  Lemma \ref{lemma-preumss}, which guarantees that the second condition in \eqref{supporting} is automatically satisfied when the first condition holds.
  \end{proof}
  \subsection{Examples}
  Take $p_0\in\mathbb H_t$ be such that $p_0\not = -p_0^\ct$, and let $D=0$,
  \[
A=\begin{pmatrix}p_0&0\\ 0&-p_0^\ct\end{pmatrix},\quad C=B^\ct=\begin{pmatrix}1&1\end{pmatrix}\quad{\rm and}\quad H=\begin{pmatrix}0&1\\1&0\end{pmatrix}.
\]
Then, equations \eqref{sna1}-\eqref{sna4} are in force, and the corresponding function $\phi$ is given by
\[
\phi(x)=x(1-xp_0)^{-1}+x(1+p_0^\ct)^{-1}
\]
and does not admit non-trivial $\ct$-unitary minimal additive decompositions since $p_0\not= -p_0^\ct$. When $p_0=x_0$ is real we have
\[
\phi(x)=2x(1-xx_0)^{-1}=x(1-xx_0)^{-1}+x(1-xx_0)^{-1}.
\]

One can obtain quite a general family of examples in the following way:

\begin{proposition}
  \label{ratio-line}
  Take $\psi(x)=D+xC(I_{\mathbb H_t^N}-xA)^{-1}$, and let
  \begin{equation}
\phi(x)=\psi(x)-\psi(-x)^\ct.
    \label{phi-psi-1}
\end{equation}
Then $\psi$ satisfies \eqref{cara-sym}. It admits a (in general not minimal) realization
\begin{equation}
  \label{psi-cal-1}
  \phi(x)=\mathcal D+x\mathcal C(I_{\mathbb H_t^{2N}}-x\mathcal A)^{-1}\mathcal B
  \end{equation}
  with
  \begin{eqnarray}
    \label{a-0111-a}
    \mathcal A&=&\begin{pmatrix}A&0\\0&-A^\ct\end{pmatrix},\\
    \mathcal B&=&\begin{pmatrix}B \\ C^\ct\end{pmatrix},\\
    \mathcal C&=&\begin{pmatrix}C&B^\ct\end{pmatrix}\\
    \mathcal D&=&D-D^\ct,
                  \label{d-0111-d}
    \end{eqnarray}
                    and satisfies \eqref{sna1}-\eqref{sna4} with
                    \begin{equation}
                      H=\begin{pmatrix} 0&I_{\mathbb H_t^N}\\ I_{\mathbb H_t^N}&0\end{pmatrix}.
                      \end{equation}
Finally $H$ is the associated symmetric matrix when the realization \eqref{a-0111-a}-\eqref{d-0111-d} is minimal.
          \end{proposition}

\begin{proof}
We have                                                                                 
  \[
  \begin{split}
    \psi(x)-\psi(-x)^\ct&=D-D^\ct+xC(I_{\mathbb H_t^{N}}-xA)^{-1} B+xB^\ct (I_{\mathbb H_t^{N}}+xA^\ct)^{-1}C^\ct\\
    &=D-D^\ct+x\begin{pmatrix}C&B^\ct\end{pmatrix}\left(I_{\mathbb H_t^{2N}}-\begin{pmatrix}A&0\\0&-A^\ct\end{pmatrix}\right)^{-1}\begin{pmatrix} B\\ C^\ct\end{pmatrix}
\end{split}
\]
which is equal to \eqref{psi-cal-1} with the realization given by \eqref{a-0111-a}-\eqref{d-0111-d}. The other claims are easily checked.
\end{proof}

            We conclude with:

            \begin{proposition}
  A $\mathbb H_t$ matrix-valued rational function satisfies \eqref{cara-sym} if and only if it can be written as
  \eqref{phi-psi-1} for some rational function $\psi$.            
\end{proposition}

          \begin{proof}
            If $\phi$ satisfies \eqref{cara-sym} then it can be written in the form \eqref{phi-psi-1} with $\psi=\frac{\phi}{2}$.
            Conversely, any function of the form
            \eqref{phi-psi-1} satisfies \eqref{cara-sym}.
            \end{proof}

\section{$\ct$-rational function $\ct$-anti-symmetric on the unit circle}
\setcounter{equation}{0}
In the present setting a matrix-valued rational function $\phi$ will be called $\ct$-anti-symmetric on the unit circle if it satisfies
\begin{equation}
  \label{phi-phi}
  \phi(x)=-(\phi(1/x))^\ct
\end{equation}
at those real points $x$ where both sides are defined. See Remark \ref{remark-r-t-2} for the motivation for this constraint.
We proceed along the lines of Section \ref{phi--phi-ct} and  reduce the study to the case of a function $J$-unitary on the unit circle
(meaning \eqref{J-TTT}) via \eqref{J-T} to study realizations and additive decompositions of such matrix-valued functions.

\subsection{Realization theorem}
\begin{theorem}
Let $\phi$ be a $\mathbb H_t^{n\times n}$-valued rational function regular at the origin and at infinity, and let
$\phi(x)=D+xC(I_{\mathbb H_t^{N\times N}}-xA)^{-1}B$ be a minimal realization of $\phi$. Then, $\phi$
satisfies \eqref{phi-phi}
at those real points $x$ where both sides are defined if and only if there exists an invertible $\ct$-symmetric matrix $H\in\mathbb H_t^{N\times N}$ such that
  \begin{eqnarray}
    \label{a-1}
    A&=&H(A^\ct)^{-1}H^{-1}\\
    D+D^\ct&=&B^\ct HB\\
    B&=&AH^{-1}C^\ct.
         \label{a-3}
  \end{eqnarray}
 The matrix $H$ is invertible and uniquely determined from the realization.
  \end{theorem}

  \begin{proof}
    Let, as in \eqref{J-T}
    \[
      \begin{split}
    J=\begin{pmatrix}0&I_{\mathbb H_t^n}\\ I_{\mathbb H_t^n}&0\end{pmatrix}\quad{\rm and}\quad
    T(x)=\begin{pmatrix}I_{\mathbb H_t^n}&\phi(x)\\0&I_{\mathbb H_t^n}\end{pmatrix}.
  \end{split}
  \]
  Condition \eqref{phi-phi} is equivalent to $T(x)JT(1/x)^\ct=J$. Applying Theorem \ref{unit-circle} to the realization \eqref{Real-T} we obtain that
  there exists a uniquely defined invertible $\ct$-symmetric matrix $H\in\mathbb H_t^{N\times N}$ such that
      \begin{equation}
        \begin{pmatrix}          A&0&B\\   C&I&D\\ 0&0&I \end{pmatrix}^\ct        \begin{pmatrix}          -H&0&0\\ 0&0&I\\0&I&0\end{pmatrix}
        \begin{pmatrix}          A&0&B\\   C&I&D\\ 0&0&I \end{pmatrix}=        \begin{pmatrix}          -H&0&0\\ 0&0&I\\0&I&0\end{pmatrix}.
 \end{equation}

 This equation can rewritten as the system
 \begin{eqnarray}
   \label{ahah123}
   A^\ct HA&=&H\\
   \label{ahah1234!}
     -A^\ct HB+C^\ct&=&0\\
   D+D^\ct&=&B^\ct HB,
              \label{ahah1234}
   \end{eqnarray}
equivalent to \eqref{a-1}-\eqref{a-3}.
\end{proof}

  We note that \eqref{ahah123} means that $A$ is unitary with respect to the form \eqref{[]H}.\\

As a corollary  of the previous theorem we have:

\begin{corollary}
In the above notations and setting, the following equalities hold:
\begin{eqnarray}
  \label{real-phi-988}
  \phi(x)&=&\frac{D-D^\ct}{2}+\frac{1}{2}C(I_{\mathbb H_t^N}+xA)(I_{\mathbb H_t^N}-xA)^{-1}H^{-1}C^\ct\\
  \frac{\phi(x)+\phi(y)^\ct}{1-xy}&=&C(I_{\mathbb H_t^N}-xA)^{-1}H^{-1}(I_{\mathbb H_t^N}-yA^\ct)^{-1}C^\ct.
                                      \label{real-phi-987}
               \end{eqnarray}
  \end{corollary}
  \begin{proof}
        Making use of \eqref{a-1}-\eqref{a-3} we can write
    \[
      \begin{split}
                \phi(x)&=\frac{D-D^\ct}{2}+\frac{D+D^\ct}{2}+xC(I_{\mathbb H_t^N}-xA)^{-1}B\\
        &=\frac{D-D^\ct}{2}+\frac{B^\ct HB}{2}+xC(I_{\mathbb H_t^N}-xA)^{-1}AH^{-1}C^\ct\\
        &=\frac{D-D^\ct}{2}+\frac{CH^{-1}A^\ct H AH^{-1}C}{2}+xC(I_{\mathbb H_t^N}-xA)^{-1}AH^{-1}C^\ct\\
        &=\frac{D-D^\ct}{2}+\frac{CH^{-1}C}{2}+xC(I_{\mathbb H_t^N}-xA)^{-1}AH^{-1}C^\ct\\
                        &=\frac{D-D^\ct}{2}+\frac{1}{2}C\left(I_{\mathbb H_t^N}+2xC(I_{\mathbb H_t^N}-xA)^{-1}A\right)H^{-1}C^\ct
        \end{split}
      \]
      from which \eqref{real-phi-988} follows. We now prove \eqref{real-phi-987}. We first note that
      \eqref{ahah123} implies that $A^\ct HA H^{-1}=I_{\mathbb H_t^N}$, and so $HAH^{-1}A^\ct=I_{\mathbb H_t^N}$, that is
        \begin{equation}
          \label{inv}
          AH^{-1}A^\ct=H^{-1}.
          \end{equation}
We now write, using \eqref{real-phi-988}:
\[
      \begin{split}
        2(\phi(x)+\phi(y)^\ct)&\\
        &\hspace{-2cm}=C(I_{\mathbb H_t^N}+xA)(I_{\mathbb H_t^N}-xA)^{-1}H^{-1}C^\ct+CH^{-1}(I_{\mathbb H_t^N}+yA^\ct)(I_{\mathbb H_t^N}-yA^\ct)^{-1}C^\ct\\
        &\hspace{-2cm}=C(I_{\mathbb H_t^N}-xA)^{-1}\left\{\spadesuit\right\}  (I_{\mathbb H_t^N}-yA^\ct)^{-1}C^\ct,
      \end{split}
    \]
    with
    \[
      \begin{split}
        \spadesuit&=(I_{\mathbb H_t^N}+xA)H^{-1}(I_{\mathbb H_t^N}-yA^\ct)+(I_{\mathbb H_t^N}-xA)H^{-1}(I_{\mathbb H_t^N}+yA^\ct)\\
        &=2\left(H^{-1}-xyAH^{-1}A^\ct\right)\\
          &=2H^{-1}(1-xy)
      \end{split}
    \]
    where we have used \eqref{inv}. The result follows.
      \end{proof}

    For such formulas in the quaternionic setting see \cite{acls_milan,zbMATH06658818}.
\subsection{$\ct$-anti-symmetric minimal additive decompositions}

For the notion of $\ct$ anti symmetric matrix see Definition \ref{anti-not-anti}.

\begin{theorem}
  Let $\phi$ satisfying \eqref{phi-phi},  with minimal realization $\phi(x)=D+xC(I_{\mathbb H_t^N}-xA)^{-1}B$ and associated $\ct$-symmetric matrix
  $H\in\mathbb H_t^{N\times N}$. Then, up to an additive $\ct$-anti-symmetric decomposition of $D-D^\ct$
  into two $\ct$-anti-symmetric matrices,
\[
  D-D^\ct=M_1+M_2, \quad where\quad M_j+M_j^\ct=0,\,\, j=1,2,
\]
there is a one-to-one correspondence between $\ct$-skew self-adjoint minimal additive
  decompositions of $\phi$ and $H$ non-degenerate $A$-invariant subspaces  of $\mathbb H_t^N$.
\end{theorem}

\begin{proof}
  As in the proof of Theorem \ref{thm-line-decomp} we have $A^\times =A$, and in the present case equation \eqref{a-1} (or, more precisely, \eqref{inv})
together with  \eqref{Adelta} means that $A$ is $H$-$\ct$
unitary,
\[
  AA^\Delta=         I_{\mathbb H_t^N},
\]
and so, by Lemma \ref{lemma-preumss}, the second condition in \eqref{supporting} is automatically satisfied when the first condition holds.
  \end{proof}
  \subsection{Examples}
  \begin{example}
Take $p_0=a+bj_t \in\mathbb H_t$ such that $p_0p_0^\ct=1$ (i.e. $|a|^2-t|b|^2=1$, see Lemma \ref{unit-p-0}), and take in
\eqref{ahah123}-\eqref{ahah1234} $A=p_0$ and $H=B=1$. Then $C=p_0$ and $D=\frac{1}{2}$ and the corresponding function $\phi$ is
\begin{equation}
  \label{p!!!p}
\frac{1}{2}+xp_0(1-xp_0)^{-1}=\frac{1}{2}(1+xp_0)(1-xp_0)^{-1}.
\end{equation}
We can add to $\phi$ any element $q_0\in\mathbb H_t$ such that $q_0+q_0^\ct=0$, i.e. (see Lemma \ref{unit-p-0-0})
of the form
\[
q_0=ix_0+j_tz_0,\quad x_0\in\mathbb R\quad{\rm and}\quad z_0\in\mathbb C.
\]
\end{example}
More generally:

\begin{example}
  Let $p_1,\ldots,p_N\in\mathbb H_t$ satisfy $p_jp_j^\ct=1$, $j=1,\ldots, N$, and set
  \[
    A={\rm diag}\,(p_1,\ldots, p_N),\quad H=I_{\mathbb H_t^N},\quad C=\begin{pmatrix}1&1&\cdots &1\end{pmatrix}\in\mathbb H_t^N.
  \]
  Then, from \eqref{ahah1234!} we get
  \[
    B=\begin{pmatrix}p_1\\ p_2\\ \vdots\\ p_N\end{pmatrix}
  \]
  and $D=\frac{N}{2}$ satisfies \eqref{ahah1234}.
  We get
  \[
    \phi(x)=\sum_{j=1}^N\left(\frac{1}{2}+x(1-xp_j)^{-1}p_j\right)
.
    \]
  \end{example}

As another example we choose
\[
H=\begin{pmatrix}0&1\\1&0\end{pmatrix}\quad{\rm and}\quad A=\begin{pmatrix}p_1&0\\0&p_1^{-\ct}\end{pmatrix},
\]
where $p_1\in\mathbb H_t$ is invertible. Then, \eqref{ahah123} is satisfied. We furthermore take 
\[
B=\begin{pmatrix}1\\1\end{pmatrix},\quad C=\begin{pmatrix}p_1&p_1^{-\ct}\end{pmatrix},\quad {\rm and}\quad D=1.
\]
Then, \eqref{ahah1234!} and \eqref{ahah1234} are also satisfied, and the corresponding $\phi$ is equal to
\[
  \begin{split}
    \phi(x)&=1+xp_1(1-xp_1)^{-1}+xp_1^{-\ct}(1-xp_1^{-\ct})^{-1}\\
    &=1+xp_1(1-xp_1)^{-1}+xp_1(p_1^\ct-x)^{-1}.
    \end{split}
\]
For $p_1p_1^\ct=1$ we get back (twice) \eqref{p!!!p}. When $p_1p_1^\ct\not=1$, the above function $\phi$ lacks non-trivial $\ct$-anti-symmetric decompositions.\\

In a way similar to Proposition \ref{ratio-line} a quite general family of examples can be obtained as follows.

\begin{proposition}
  Take $\psi(x)=D+xC(I_{\mathbb H_t^N}-xA)^{-1}$ with $A$ invertible, and let
  \begin{equation}
\phi(x)=\psi(x)-\psi(1/x)^\ct.
    \label{phi-psi}
\end{equation}
Then $\psi$ satisfies \eqref{phi-phi}. It admits a (in general not minimal) realization
\begin{equation}
  \label{psi-cal}
  \phi(x)=\mathcal D+x\mathcal C(I_{\mathbb H_t^{2N}}-x\mathcal A)^{-1}\mathcal B
  \end{equation}
  with
  \begin{eqnarray}
    \label{a-0111}
    \mathcal A&=&\begin{pmatrix}A&0\\0&(A^\ct)^{-1}\end{pmatrix},\\
    \mathcal B&=&\begin{pmatrix}B \\ (A^\ct)^{-1}C^\ct\end{pmatrix},\\
    \mathcal C&=&\begin{pmatrix}C&B^\ct(A^\ct)^{-1}\end{pmatrix},\\
    \mathcal D&=&D-D^\ct+B^\ct(A^\ct)^{-1}C^\ct,
                  \label{d-0111}
    \end{eqnarray}
                    and satisfies \eqref{ahah123}-\eqref{ahah1234} with
                    \begin{equation}
                      H=\begin{pmatrix} 0&I_{\mathbb H_t^N}\\ I_{\mathbb H_t^N}&0\end{pmatrix}.
                      \end{equation}
Finally $H$ is the associated $\ct$-symmetric matrix when the realization \eqref{a-0111}-\eqref{d-0111} is minimal.
          \end{proposition}

\begin{proof}
We have                                                                                 
  \[
  \begin{split}
    \psi(x)-\psi(1/x)^\ct&=D-D^\ct+xC(I_{\mathbb H_t^{2N}}-xA)^{-1} B-B^\ct (xI_{\mathbb H_t^{2N}}- A^\ct)^{-1}C^\ct\\
    &=D-D^\ct+xC(I_{\mathbb H_t^{2N}}-xA)^{-1} B+B^\ct (A^{\ct})^{-1} (I_{\mathbb H_t^{N}}- x(A^\ct)^{-1})^{-1}C^\ct\\
    &=D-D^\ct+B^\ct (A^\ct)^{-1}C^\ct+xC(I_{\mathbb H_t^{N}}-xA)^{-1} B+\\
    &\hspace{5mm}+B^\ct (A^{\ct})^{-1} (I_{\mathbb H_t^{N}}- x(A^\ct)^{-1})^{-1}C^\ct-B^\ct(A^\ct)^{-1}C^\ct\\
    &=D-D^\ct+B^\ct (A^\ct)^{-1}C^\ct+xC(I_{\mathbb H_t^{N}}-xA)^{-1} B+\\
    &\hspace{5mm}+xB^\ct (A^{\ct})^{-1} (I_{\mathbb H_t^{N}}- x(A^\ct)^{-1})^{-1}(A^{\ct})^{-1}C^\ct\\
    &=D-D^\ct+B^\ct (A^\ct)^{-1}C^\ct+\\
    &\hspace{5mm}+x\begin{pmatrix}C &B^{\ct}(A^{\ct})^{-1}\end{pmatrix}\left(I_{\mathbb H_t^{2N}}-x\begin{pmatrix}A&0\\0&(A^\ct)^{-1}\end{pmatrix}\right)^{-1}\begin{pmatrix} B\\ (A^{\ct})^{-1}C^\ct\end{pmatrix},
\end{split}
\]
which is equal to \eqref{psi-cal} with the realization given by \eqref{a-0111}-\eqref{d-0111}. The other claims are readily verified.
\end{proof}

\begin{example}
Taking in the previous proposition $\psi(x)=(1-xp_0)^{-n}$ we get
\[
\varphi(x)=(1-xp_0)^{-n}-x^n(x-p_0^\ct)^{-n}.
\]
\end{example}

We conclude with:

            \begin{proposition}
  A $\mathbb H_t$ matrix-valued rational function satisfies \eqref{phi-phi} if and only if it can be written as
  \eqref{phi-psi} for some rational function $\psi$.            
\end{proposition}

          \begin{proof}
       If $\phi$ satisfies \eqref{phi-phi} then it can be written in the form \eqref{phi-psi} with $\psi=\frac{\phi}{2}$. Conversely, any function of the form
            \eqref{phi-psi} satisfies \eqref{phi-phi}.
            \end{proof}
\bibliographystyle{plain}
\bibliography{all}
  \end{document}